\documentclass[a4paper,leqno]{siamltex}

\usepackage[thinspace,amssymb,thinqspace,thinspace]{SIunits}
\usepackage{amsmath}
\usepackage{stmaryrd}
\usepackage{amssymb}
\usepackage{graphicx}
\usepackage{subfigure}
\usepackage{booktabs}
\usepackage{xspace}

\usepackage{cite}
\usepackage{soul}
\usepackage{multirow}
\usepackage{algorithm}
\usepackage[]{algorithmicx}
\usepackage[]{algpseudocode}
\algrenewcommand\algorithmicrequire{\textbf{Input:}}
\algrenewcommand\algorithmicensure{\textbf{Output:}}
\algrenewcommand\algorithmicforall{\textbf{for each}}

\usepackage[left=3cm,right=2cm,top=2cm,bottom=2cm]{geometry}
\usepackage[utf8]{inputenc}

\usepackage{tikz}
\usetikzlibrary{tikzmark,calc,decorations.pathreplacing,arrows.meta,calligraphy}
\usepackage[citecolor=black]{hyperref}

\DeclareMathOperator*{\argmin}{arg\,min}

\newcommand{\REMOVE}[1]{}

\newcommand{\StencilTwoD}[4]{ %
  \left[ %
    \begin{array}{ccccc}%
      #4 &  & #2 & & #4\\%
      &&&&\\%
      #3 &  & #1 & & #3\\%
      &&&&\\%
      #4 &  & #2 & & #4\\%
    \end{array}\right]}

\newcommand{\phantomLeft}{\left . }
\newcommand{\phantomRight}{\color{white} \right \} }

\newcommand{\gb}{\mathbf{g}}
\newcommand{\eb}{\mathbf{e}}
\newcommand{\vb}{\mathbf{v}}

\newcommand{\Ep}{E^\prime}
\newcommand{\Vp}{V^\prime}
\newcommand{\gbp}{\gb^\prime}

\newcommand{\norm}[1]{\left\lVert#1\right\rVert}
\newcommand{\abs}[1]{\left\lvert#1\right\rvert}
\newcommand{\neuralnet}{ {\it N \hskip -.03in N} }
\newcommand{\subscriptneuralnet}[1]{ {\it {N \hskip -.03in N}_{\hskip -.025in {#1}}}}
\newcommand{\boldterminates}[2]{\mathbf{#1}_{*#2}^\prime}

\definecolor{darkgreen}{rgb}{0,0.4,0}
\definecolor{darkpurple}{rgb}{0.5,0,0.5}
\definecolor{clemson}{RGB}{245,102,0}

\graphicspath{{fig/}{eps/}}

\begin{document}

\title{Graph Neural Networks and Applied Linear Algebra\thanks{ 
      This work was supported by the U.S.~Department of Energy, Office of Science, Office of Advanced Scientific
       Computing Research, Applied Mathematics program, Exascale Computing Project, Early Career Research Program,
       and the SEA-CROGS project in the MMICCs program.
       Sandia National Laboratories is a multimission laboratory managed and operated by National Technology
       and Engineering Solutions of Sandia, LLC., a wholly owned subsidiary of Honeywell International, Inc.,
       for the U.S. Department of Energy's National Nuclear Security Administration under grant~DE-NA-0003525.
       This paper describes objective technical results and analysis.  Any subjective views or opinions that
       might be expressed in the paper do not necessarily represent the views of the U.S. Department of Energy
       or the United States Government. Sandia release number: SAND2023-10755O.}}
\author{
  Nicholas S. Moore\thanks{West Texas A\&M University, 2501 4th Ave, Canyon, TX 79016
    ({nmoore@wtamu.edu}).}
  \and
  Eric C. Cyr\thanks{Sandia National Laboratories, P.O. Box 5800, MS 1320, Albuquerque, NM 87185
    ({eccyr@sandia.gov}).}
  \and
  Peter Ohm\thanks{RIKEN Center for Computational Science, 7-1-26 Minatojima-minami-machi, Chuo-ku, Kobe, Hyogo, 650-0047, Japan
    ({peter.ohm@riken.jp}).}
  \and
  Christopher M. Siefert\thanks{Sandia National Laboratories, P.O. Box 5800, MS 1320, Albuquerque, NM 87185 ({csiefer@sandia.gov}).}
  \and
  Raymond S. Tuminaro\thanks{Sandia National Laboratories, P.O. Box 5800, MS 1320, Albuquerque, NM 87185 ({rstumin@sandia.gov}).}
}

\maketitle

\begin{abstract}
Sparse matrix computations are ubiquitous in scientific computing. Given the recent interest in scientific machine learning, it is natural to ask how sparse matrix computations can leverage neural networks (NN). Unfortunately, multi-layer perceptron (MLP) neural networks are typically not natural for either graph or sparse matrix computations.  The issue lies with the fact that MLPs require fixed-sized inputs while scientific applications generally generate sparse matrices with arbitrary dimensions and a wide range of different nonzero patterns (or  matrix graph vertex interconnections).  While convolutional NNs could possibly address matrix graphs where all vertices have the same number of nearest neighbors, a more general approach is needed for arbitrary sparse matrices, e.g.  arising from discretized partial differential equations on unstructured meshes.  Graph neural networks (GNNs) are one such  approach suitable to sparse matrices. The key idea is to define aggregation functions (e.g., summations) that operate on variable size input data to produce data of a fixed output size so that MLPs can be applied.  The goal of this paper is to provide an introduction to GNNs for a numerical
linear algebra audience.  Concrete GNN examples are provided to illustrate how many common linear algebra tasks can be accomplished using GNNs.  We focus on iterative and multigrid methods that employ computational kernels such as matrix-vector products, interpolation, relaxation methods, and strength-of-connection measures. Our GNN examples include cases where parameters are determined \emph{a-priori} as well as cases where parameters must be learned.  The intent with this article is to help computational scientists understand how GNNs can be used to adapt machine learning concepts to computational tasks associated with sparse matrices. It is hoped that this understanding will further stimulate data-driven extensions of classical sparse linear algebra tasks.

\end{abstract}

\section{Introduction}
Artificial intelligence (AI) and machine learning (ML) have drawn a great deal
of media attention 
 --- deep fake audio \cite{Khanjani2023},
transformer-based AI chatbots \cite{huggingfacetransformers}, stable diffusion AI
image generation and even AI Elvis \cite{AIElvis} singing a
song have all inserted themselves into popular culture.  In 
scientific realms, medical image identification (e.g. detection of cancer) has
also generated media attention exposing the potential of ML technologies.  But in the
realm of applied mathematics the impact of AI/ML has been far less
visible to the general public. The goal of this paper is to provide an introduction to graph neural networks (GNNs) and
to show how this specific  class of AI/ML algorithms can be used to represent 
(and enhance) traditional algorithms in numerical linear algebra. Associated code which implements the example GNNs is provided at \url{https://github.com/sandialabs/gnn-applied-linear-algebra/}.

Introductions to neural network models in ML and AI typically focus
on deep neural networks (DNNs) or convolutional neural networks (CNNs).
However, GNNs are can be notably more appropriate 
for many computational science tasks~\cite{shukla2022scalable}.  DNNs and CNNs generally
assume \textit{structured} input data --- a vector that is a
 fixed size, or an image where pixels are aligned along Cartesian directions. 
While structured applications do exist in
computational science (uniformly meshed problems for example), many other applications are \textit{unstructured}, 
again often associated with meshes.
Consider a drawing of an object that an engineer wishes to
simulate on a computer.  After pre-processing, the geometry of the object will be subdivided a 
mesh. 
The object may have holes, protrusions or a
disparity of features scales such that a structured, Cartesian mesh
cannot be generated. As a result, an unstructured mesh is constructed.
Application of traditional DNN and CNN tools would likely not be possible for this mesh due to the 
the unstructured nature of the data.  
However, GNN models could still be employed on this geometry for associated AI/ML tasks. 

So what is a GNN?  As suggested by the name, graph neural networks are
built on the concept of associating data with edges and vertices of a
\textit{graph}. For the above engineered object, the mesh itself can
define the graph.  
Graphs can also represent
structured data ---  for an image, pixels can be vertices and
edges can be used to represent 
neighboring pixels.  
The general nature
of unstructured graphs make them highly appropriate for a wide range
of science and engineering applications.   Beyond the meshing example,
 another example of unstructured data arises from social network graphs where
each edge represents a connection between two people. In these type of
graph networks, eigenvector 
calculations provide useful information for determining the influence
that a vertex (e.g., a person) has on the rest of the network. 
GNNs were first proposed in~\cite{gori_new_2005}
and~\cite{scarselli_graph_2009} as a means of adapting a convolutional
neural
network to graph problems.  

The issue with basic neural networks (e.g. DNNs and CNNs) that GNNs
address is that they require a fixed input/feature size. 
That is, a basic neural network can be viewed as a type of function
approximation where the number of input/feature values 
to the function is always the same when either training the network or when using the network for inference.
If a graph has a simple repeatable interaction pattern, then 
a basic neural network can be effectively applied to a fixed window (i.e., a
subset of the graph), which can be moved to address different portions
of the network.  However, general graphs have no 
such repeatable pattern.  The key difficulty lies in the fact that the number of edges adjacent to each vertex can 
vary significantly throughout the network. In the case of a social
network graph, for instance, a popular person has many 
friends while a loner might have few interactions.
To address this variability, GNNs include the notion of general aggregation functions that allow for a variable 
number of inputs. Simple aggregation function examples are summation or maximum which are well-defined
regardless of the number of inputs. Aggregation functions are used to combine information from edges
adjacent to a vertex. These functions can include a fixed number of learnable parameters. For example, 
an aggregation function might include both a summation and a maximum function 
whose results are combined in a weighted fashion using a learnable
weight. These are then combined with transformation functions
associated with edges and vertices, which can also include learnable parameters.

It is well known that many numerical linear algebra algorithms also
possess a graph structure. For instance, in Section~\ref{sec:algo_matvec} we detail how 
sparse matrix-vector multiplication transforms matrix entries as graph edge attributes and vector values as graph vertex
attributes, into a new set of graph vertex attributes representing the product. This algorithm can be written as a GNN with
precise choices of aggregation and transformation functions. We
provide additional descriptions of standard iterative methods and
components of algebraic multigrid (AMG) preconditioners
that can be re-formulated using GNNs.  Some of the considered multigrid components include matrix-vector
products, Jacobi relaxation, AMG strength-of-connection measures, and AMG interpolation operators. While the associated
GNNs do not contain trainable parameters, they illustrate the flexibility of GNNs in representing traditional
numerical linear algebra components.
The potential for advances in AMG using GNN architectures is made apparent by
including trainable parameters in the transformations. These parameters will enable learning complex nonlinear
relationships between the feature spaces, implicitly encoding the topological and numerical
properties of a sparse matrix. We explore some prototypical use cases in Section~\ref{sec:iter-methods}.
We note that more complicated
applications of GNNs for multigrid have been considered in the
literature and encourage the interested reader to look more deeply
\cite{Moore2021, luz2020,taghibakhshi2023_MGGNN,taghibakhshi2022_OptimizationBasedAMG}. As this paper is educational in nature, we are not
proposing new AI/ML powered methods that significantly outperform
existing methods, but rather demonstrating how GNNs can allow us to
modify existing numerical methods while simultaneously highlighting some of the challenges that must be
addressed.

\subsection{Developing Intuition: A Viral Example}

To make the discussion 
more accessible, we present
an example based on the spread of disease (a subject that we are all
unfortunately familiar with). The intent is to help the reader develop
intuition about 
components of the GNN algorithms
and data structure. In the text that follows we call out explicitly where we
use the example. 
If the notation and ideas are clear to the
reader, these \emph{Viral Interludes} can be skipped.  

Here, we briefly describe the setup of the viral example. We consider
a community where each individual interacts on a weekly basis with
 a fixed set of community members. 
 The one-week period over which the set of interactions occurs is 
 referred to as a cycle. Initially, within the community each individual
 has a probability of being infected by the disease. Simultaneously, a
 cure, which also (conveniently) diffuses through individual interactions,
has been distributed randomly through the population\footnote{Perhaps additional individuals seek treatment after learning about it through peer
interactions.}.
The specific rate
 of diffusion for the disease and cure, as well as the effectiveness of the
 intervention is not known in advance, and must be estimated based on observations.
The goal then is to develop a model to
assess each individual's risk of carrying the infection or carrying
the cure after a particular number of cycles using the information
about the interactions between individuals.

\section{Graph Neural Network Background}
Graph neural networks (GNN) were introduced in~\cite{gori_new_2005,scarselli_graph_2009} to overcome limitations in
convolutional neural networks (CNN) that assume a regular structure on the input data to the network. A CNN takes, as
input data, structured data that is topologically like a Cartesian mesh, and produces 
output data 
with a similarly regular structure.  Bitmap images are the usual exemplar for
CNN applications.
By contrast, the topology of
input/output
data elements for a GNN is associated with a graph. A GNN takes attributes for each vertex and for each edge
and produces a new set of vertex and edge attributes. Additionally, a set of global attributes associated with the 
entire graph can also be transformed to a new set of global attributes.

Formalizing this description, define a directed graph $\mathcal{G} =
(\mathcal{V},\mathcal{E})$, where $\mathcal{V}$ is the set of $N^v$
vertices, and $\mathcal{E}$ is the set of $N^e$ directed edges. A directed edge
is defined as an ordered pair of vertices, $(v_i,v_j) = e_{ij} \in
\mathcal{E}$ where the edge emanates from $v_i \in \mathcal{V}$ and terminates at $v_j \in \mathcal{V}$ . The set of
features/attributes (data defined on vertices are denoted by $V = \{\vb_k \in \mathbb{R}^{n_v} : k=1\ldots
N^v \}$, on edges by $E = \{\eb_{ij}  \in \mathbb{R}^{n_e} : e_{ij} \in \mathcal{E}\}$, and for the graph are $\gb \in \mathbb{R}^{n_g}$. 
While the terms ``features'' and ``attributes'' are both used
to describe data associated with edges or vertices in the literature,
following Battaglia et al. \cite{battaglia_relational_2018}, we will
prefer the term attributes. 
Table~\ref{tab:notation} summarizes notation used within this paper (including symbols that are introduced shortly).

\begin{table}[htpb]
  \centering
  \caption{Notation summary}
  \label{tab:notation}
  \begin{tabular}{cl}
    Symbol & Meaning \\
    \hline
    $\mathcal{G}$         & Graph \\
    $\mathcal{V}$         & Set of vertices in a graph \\
    $\mathcal{E}$         & Set of edges in a graph \\
    $\mathcal{N}$         & Graph neural network \\
    $\Theta$              & Graph neural network parameters \\
    $N^{v}$               & Number of vertices in a graph: $\abs{\mathcal{V}}$  \\
    $N^{e}$               & Number of edges in a graph: $\abs{\mathcal{E}}$ \\
    $V$                   & Set of input attributes defined on vertices \\
    $E$                   & Set of input attributes defined on edges \\
    $\mathbf{g}$          & Set of input global attributes on a graph \\
    $n_{v}$               & Number of input vertex attributes per vertex \\
    $n_{e}$               & Number of input edge attributes per edge \\
    $n_{g}$               & Number of input global attributes \\
    $v_{k}$               & $k^{th}$ vertex in the graph \\
    $e_{ij}$              & Directed edge that emanates from $v_i$ and terminates at $v_j$ \\
    $\mathbf{v}_k$        & Attributes associated with vertex $v_{k}$ \\
    $\textbf{e}_{ij}$     & Attributes associated with edge $e_{ij} $  \\
    $\overline{\eb_j}$    & Aggregated attributes from all edges terminating at vertex $v_{j}$ \\ 
    $\mathbb{R}^{\mathrm{var}(v)}$ & Variable input space for variadic vertex-based functions \\
    $\rho_{e \rightarrow v}()$ & Function whose inputs (outputs) are edge (vertex) attributes of a vertex \\
    $\rho_{e \rightarrow g}()$ & Function whose inputs (outputs) are edge (global) attributes of the graph \\
    $\rho_{v \rightarrow g}()$ & Function whose inputs (outputs) are vertex (global) attributes of the graph \\
    $n_{e \rightarrow v}$ & Number of output vertex attributes produced by $\rho_{e \rightarrow v}()$ \\
    $n_{e \rightarrow g}$ & Number of output global attributes produced by $\rho_{e \rightarrow g}()$ \\
    $n_{v \rightarrow g}$ & Number of output global attributes produced by $\rho_{v \rightarrow g}()$ \\
  \end{tabular}
\end{table}

\paragraph{Viral Interlude} For our viral example each vertex in the graph $v_i \in\mathcal{V}$ represents the $i^{th}$ individual in the community. The initial vertex attribute, prior to the first cycle, is the 
vector $\mathbf{v}_i \in [0,1]^2 \in\mathbb{R}^{n_v}$ where $n_v=2$. The first component is the probability the individual contains the
cure, while the second component is the probability the individual contains the disease. Similarly, the existence of an edge $e_{ij} \in \mathcal{E}$ 
indicates an interaction between individual $i$ and $j$. Initially, this edge contains as attributes the length of the interaction, which is important in determining both the spread of the
infection and the cure, thus $\mathbf{e}_{ij} \in \mathbb{R}^{n_e}$ where $n_e=1$. 

\subsection{Graph Neural Network Function}

A GNN, denoted by $\mathcal{N}$, is a parameterized function that acts on the attributes of a graph and produces new attributes while
maintaining topology:
\begin{equation}
\mathcal{N}\left(\mathcal{G}, (V, E, \gb); \Theta\right) = \left \{ V', E', \gb' \right \} ,
\end{equation}
where $\Theta$ are the GNN parameters, and a tick denotes the output
attributes. Note that the sizes of the attributes can change as a
result of the neural network action; for example $n_v, n_e, n_g \mapsto n_{v^\prime}, n_{e^\prime},
n_{g^{\prime}}$ where $n_{v}$ can be different from $n_{v^{\prime}}$, $n_{e}$ can be different from
$n_{e^{\prime}}$, and $n_{g}$ can be different from $n_{g^{\prime}}$.
This feature of GNNs will be used in the Viral Interludes, though not
in the examples in Section~\ref{sec:iter-methods}.
Even if the number of attributes change after the
GNN, the topology of the
graph where attributes reside is not perturbed by the graph neural network. However,
edges and vertices can be marked as on or off in the attributes, though the fundamental 
topology doesn't change. Also note that the GNN architectures discussed here apply only to
directed graphs. Undirected graphs can be represented as well by treating an undirected edge
as a pair of directed edges, one in each direction. 

In this paper, we focus on GNN models based on message passing that were originally proposed
in~\cite{gilmer_neural_2017}. This model decomposes the action of the neural network into a sequence
of layers, that themselves are parameterized functions of the graph topology and attributes. For
instance, a two layer GNN with layers denoted as $\mathcal{L}_{1}$ and $\mathcal{L}_{2}$ is
\begin{equation}
\mathcal{N}(\mathcal{G}, (V, E, \gb); \{\Theta_1,\Theta_2\})
   :=                \mathcal{L}_{2}(\mathcal{G},
                   \mathcal{L}_{1}(\mathcal{G},(V,E,\gb);\Theta_1);
                \Theta_2).
\end{equation}
The specifics of layer architectures are discussed in Section~\ref{sec:mpl}.

\paragraph{Viral Interlude} 
For our viral model, a GNN function $\mathcal{N}$ represents the
evolution of the disease within this community from the initial state
to a final state. $\mathcal{N}$ is constructed through composition of
the layers, where a layer $\mathcal{L}_i$ describes an evolution over
a single cycle. In a cycle new attributes for each individual (the set
$V'$) are computed by considering a combination of the vertex
attributes $V$, and the time of interaction specified by the edge
attribute $E$. The parameters $\Theta$ are
not known \textit{a priori} and must be calibrated based on empirical data about the evolution of the
outbreak. (see Section~\ref{sec:GNN-training}). In a GNN, the parameter $\Theta$ encompasses these model unknowns.
The layer may also require edge
updates based on changing interaction times. Further, if the number of values
representing an edge or vertex state changes as a result of applying a layer
then we have $n_v' \neq n_v$ and/or $n_e'\neq n_e$. For instance, an
additional attribute could be added by the GNN layer for each vertex
specifying the survival (or not) of an individual. This boolean
attribute could have the effect of virtually removing an individual
from the graph in subsequent layers of the GNN
function. 

\subsection{Message Passing Layer} \label{sec:mpl}

While several different variations of message passing GNNs exist, we will describe and utilize the GNN from \cite{battaglia_relational_2018}. The parameterization of a GNN layer is
defined by three aggregation functions and three update functions. Because graphs can have varying
topology, a method for combining attributes from multiple entities of the same type is necessary.
Formally, these
aggregation functions are
\begin{equation}
	\rho_{e \rightarrow v}(\boldterminates{e}{k}): \mathbb{R}^{\mathrm{var}(v)} \rightarrow \mathbb{R}^{n_{e \rightarrow v}}, \quad
	\rho_{e \rightarrow g}(E^{\prime} ): \mathbb{R}^{N^{e} \cdot n_e} \rightarrow \mathbb{R}^{n_{e \rightarrow g}}, \quad
	\rho_{v \rightarrow g}(V^{\prime} ): \mathbb{R}^{N^{v} \cdot n_v} \rightarrow \mathbb{R}^{n_{v \rightarrow g}}
\end{equation}
where $\boldterminates{e}{k}$ denotes the set of updated
attributes associated with all edges that terminate at vertex
$v_k$. 
The function $\rho_{e \rightarrow v}()$ is applied at each vertex (i.e., for each $ k = 1 , \dots, N_{v} $)
and the notation $\mathbb{R}^{var(v)}$ denotes that the $\rho_{e \rightarrow v}$ function is
variadic: the input space is a finite (but variable) number of edge attribute vectors. This function
takes the set of attributes from all edges which terminate at a single vertex and aggregates them
together into $n_{e \rightarrow v}$ attributes that are collected at this vertex. The other two aggregation functions, $\rho_{e
\rightarrow g}$ and $\rho_{v \rightarrow g}$, are similar but combine updated attributes 
$E^\prime$ associated with all graph edges and updated attributes $V^\prime$
associated with all graph vertices respectively.
Some simple examples of aggregation functions include
summation, minimum, and maximum. The aggregated attributes are then used to update the attributes of
the vertices and the graph.

The update functions are
\begin{align}
\phi_e(\eb_{ij}, \vb_{i}, \vb_{j}, \gb)&: \mathbb{R}^{n_e} \times \mathbb{R}^{n_v}\times  \mathbb{R}^{n_v}\times  \mathbb{R}^{n_g} \rightarrow \mathbb{R}^{n_e'}, \\
\phi_v(\vb_{i}, \overline{\eb}_i, \gb)&: \mathbb{R}^{n_v} \times  \mathbb{R}^{n_{e \rightarrow v}} \times \mathbb{R}^{n_g}\rightarrow \mathbb{R}^{n_v'}, \\
\phi_g(\gb, \rho_{e \rightarrow g}(E^{\prime}), \rho_{v \rightarrow g}(V^{\prime} ))  &: \mathbb{R}^{n_g} \times \mathbb{R}^{n_{e \rightarrow g}} \times  \mathbb{R}^{n_{v \rightarrow g}} \rightarrow \mathbb{R}^{n_g'}
\end{align}
where $\overline{\eb}_i$ are edge attributes that have been aggregated to vertex $i$ (i.e., $\overline{\eb}_i = \rho_{e \rightarrow v} (\boldterminates{e}{i})$). 
Each of these functions (with any implied parameters), takes, as input, the attributes of a graph
entity and the attributes of entities in its neighborhood and transforms them into an updated
attribute associated with the original graph entity. 
For instance $\phi_e$ is associated with
updating an edge's attributes in the graph. The edge's neighborhood includes two vertices, and the graph itself. 
Therefore, it takes as input the attributes of the edge, its two
neighboring vertices and the graph's attributes. The vertex update, $\phi_v$, is associated with a
vertex in the graph. The neighborhood of the vertex includes the graph itself and all
the edges which terminate at the vertex. The
number of terminating edges can vary from vertex to vertex, which necessitates use of the variadic
aggregation function $\rho_{e \rightarrow v}$ that combines all the edge attributes together. 
The
original attributes of the vertex, the aggregated edge attributes, and the global attributes are
then used as arguments to the vertex update function. A similar process occurs in the global update
$\phi_g$ as well, where all updated edge attributes are aggregated
together and all updated vertex attributes are aggregated together and used
to update the global graph attributes.

With the the aggregation and update functions described, we can now explain how they are combined in
Algorithm~\ref{algo:GNBlock} 
\begin{algorithm}
\caption{Computation of a Graph Network Layer}\label{algo:GNBlock}
\begin{algorithmic}[0]
\Require Graph $G$, vertex attributes $V = \{\vb_{j} : j = 1, \ldots,
N^{v}\}$, edge attributes $E = \{\eb_{ij} : e_{ij} \in \mathcal{E}\}$
and global attributes $\gb$
\Ensure updated edge attributes $E^\prime$, vertex attributes
$V^\prime$ and global attributes $\gb^\prime$

\\\hrulefill

  \Function{$(\Vp,\Ep,\gbp)$ = GraphNetworkLayer}{$G,V,E,\gb$}

  \For{$\eb_{ij} \in E $} \Comment{For each edge}
    \State $\eb_{ij}^\prime \gets \phi_e(\eb_{ij}, \vb_i, \vb_j, \gb)$ \Comment{Update edge features} 
  \EndFor
  \For{$k \in \{1\ldots N^v\}$} \Comment{For each vertex}
  \State $\overline{\eb}_k^\prime \gets \rho_{e\rightarrow v}\left(\boldterminates{e}{k}\right)$
  \Comment{Aggregate all edge features for edges terminating at $v_k$}
  \State $\vb_k^\prime \gets \phi_{v}(\vb_k, \overline{\eb}_k^\prime, \gb)$
  \Comment{Update vertex features}
  \EndFor
  \State Let $V^\prime = \{\vb_k^\prime\}_{k=1:N^v}$ 
  \State Let $E^\prime = \{\eb_{ij}^\prime\}_{e_{ij} \in \mathcal{E}}$ 
  \State $\gb^\prime \gets \phi_{g}(\gb, ~\rho_{e \rightarrow g}(E^\prime),~\rho_{v \rightarrow g}(V^\prime))$  \Comment{Update global features}
  \State \Return $(E^\prime, V^\prime, \gb^\prime)$
  \EndFunction
\end{algorithmic}
\end{algorithm}

to compute the action of a GNN layer. The first {\bf for} loop in the
algorithm transforms input edge attributes using the update function $\phi_e$. This is depicted graphically in
Figure~\ref{fig:edge-update} with the initial state of the graph shown on the left
graph and the final state shown on the right. For a single edge $e_{12}$, its attributes, along with those of the
neighboring vertices, and the graph are the input to $\phi_e$ yielding an updated edge
attribute $\eb_{12}^\prime$. 
While only one update is shown, all edges are updated (potentially in parallel) using the
same update function.
\definecolor{color1}{RGB}{135,54,0}  
\definecolor{color2}{RGB}{81,148,37} 
\definecolor{color3}{RGB}{0,169,212} 
\definecolor{color4}{RGB}{65,47,186} 
\definecolor{color5}{RGB}{199,50,90} 
\begin{figure}
  \caption{Edge feature update for $\eb_{12} $. The concatenated attributes for the edge, neighboring vertices, and
    global graph are input to the vertex update function $\phi_e$. The output $\eb_{12} ^\prime$ corresponds to 
    updated edge attributes.  This process (not shown) is repeated for all edges in parallel, using the same
    update function. }\label{fig:edge-update}
\vskip .25in
  \centering
  \begin{tikzpicture}[->,auto,
    thick,main node/.style={circle,draw,font=\sffamily\normalsize}, scale=.6]

    \node[color1, main node] (1) at (-2,2)   {$\vb_1$};
    \node[color2, main node] (2) at (2, 2)   {$\vb_2$};
    \node[main node] (3) at (-2, -2) {$\vb_3$};
    \node[main node] (4) at (2, -2)  {$\vb_4$};

    \path[every node/.style={font=\sffamily\small}]
    (1) edge[color3] (2)
    (1) edge[black] (3)
    (1) edge[black, bend left] (4)
    (2) edge[black] (4)
    (3) edge[black] (4)
    (4) edge[black, bend left] (1);

    \node[color3, scale=1.0] at (0, 2.5) {$\eb_{12} $};
    \node[scale=1.0] at (-2.5, 0) {$\eb_{13}$};
    \node[scale=1.0] at (2.6, 0) {$\eb_{24} $};
    \node[scale=1.0] at (0, -2.5) {$\eb_{34}$};
    \node[scale=1.0] at (.4, .4) {$\eb_{14} $};
    \node[scale=1.0] at (-.4, -.4) {$\eb_{41} $};

    \node[scale=1.0] at (4.5, 6) {$\displaystyle \phi_{e}\left( \begin{bmatrix} \textcolor{color3}{\eb_{12} } \\ \textcolor{color1}{\vb_1} \\ \textcolor{color2}{\vb_2} \\ \gb \end{bmatrix} \right) = \textcolor{color4}{\eb_{12} ^\prime}$};

    \node[main node] (5) at (7, 2)   {$\vb_1$};
    \node[main node] (6) at (11, 2)   {$\vb_2$};
    \node[main node] (7) at (7, -2) {$\vb_3$};
    \node[main node] (8) at (11, -2)  {$\vb_4$};

    \path[every node/.style={font=\sffamily\small}]
    (5) edge[color4] (6)
    (5) edge[black] (7)
    (5) edge[black, bend left] (8)
    (6) edge[black] (8)
    (7) edge[black] (8)
    (8) edge[black, bend left] (5);

    \node[color4, scale=1.0] at (9, 2.5) {$\eb_{12}^\prime$};
    \node[scale=1.0] at (6.5, 0) {$\eb_{13} $};
    \node[scale=1.0] at (11.6, 0) {$\eb_{24} $};
    \node[scale=1.0] at (9, -2.5) {$\eb_{34} $};
    \node[scale=1.0] at (9.4, .4) {$\eb_{14} $};
    \node[scale=1.0] at (8.6, -.4) {$\eb_{41} $};

    \draw[-{Triangle[width=18pt, length=8pt]}, line width=18pt] (3.5,0) -- (5.6,0);

    \node[white, scale=1.0] at (4.45, 0) {$\phi_{e}$ };
  \end{tikzpicture}
\vskip -.05in
\end{figure}
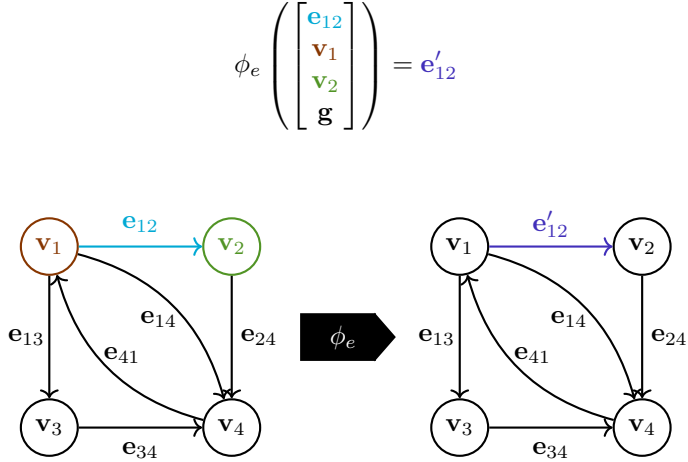

\begin{figure}
  \caption{Vertex feature update for $\vb_4$. The updated attributes for incoming edges are
    aggregated using $\rho_{e\rightarrow v}$ yielding aggregated edge information
    $\mathbf{\overline{e}}_4$. The aggregated attributes are concatenated with vertex
    attributes and global attributes to define the input to the vertex update function $\phi_v$. The
    output $\vb_4^\prime$ corresponds to updated vertex attributes. 
    This process is repeated for all vertices in parallel (not shown), using the same aggregation
    and update functions. }
    \label{fig:vertex-update}
\vskip .1in
  \centering
  \begin{tikzpicture}[->,auto,
    thick,main node/.style={circle,draw,font=\sffamily\normalsize}, scale=.6]
    \node[main node] (1) at (-11,2)   {$\vb_1$};
    \node[main node] (2) at (-7, 2)   {$\vb_2$};
    \node[main node] (3) at (-11, -2) {$\vb_3$};
    \node[teal, main node] (4) at (-7, -2)  {$\vb_4$};

    \path[every node/.style={font=\sffamily\small}]
    (1) edge[black] (2)
    (1) edge[black] (3)
    (1) edge[color3, bend left] (4)
    (2) edge[color1] (4)
    (3) edge[color2] (4)
    (4) edge[black, bend left] (1);

    \node[scale=1.0] at (-9, 2.5) {$\eb_{12} ^\prime$};
    \node[scale=1.0] at (-2.5-9, 0) {$\eb_{13} ^\prime$};
    \node[color1, scale=1.0] at (2.6-9, 0) {$\eb_{24} ^\prime$};
    \node[color2, scale=1.0] at (-9, -2.5) {$\eb_{34} ^\prime$};
    \node[color3, scale=1.0] at (0.4-9, .4) {$\eb_{14} ^\prime$};
    \node[scale=1.0] at (-9.4, -.4) {$\eb_{41} ^\prime$};

    \node[main node] (1) at (-2,2)   {$\vb_1$};
    \node[main node] (2) at (2, 2)   {$\vb_2$};
    \node[main node] (3) at (-2, -2) {$\vb_3$};
    \node[main node] (4) at (2, -2)  {$\textcolor{color5}{\overline{\eb}_{4}}, \textcolor{teal}{\vb_4}$};

    \path[every node/.style={font=\sffamily\small}]
    (1) edge[black] (2)
    (1) edge[black] (3)
    (1) edge[black, bend left] (4)
    (2) edge[black] (4)
    (3) edge[black] (4)
    (4) edge[black, bend left] (1);

    \node[scale=1.0] at (0, 2.5) {$\eb_{12}^\prime$};
    \node[scale=1.0] at (-2.5, 0) {$\eb_{13}^\prime$};
    \node[scale=1.0] at (2.6, 0) {$\eb_{24}^\prime$};
    \node[scale=1.0] at (0, -2.5) {$\eb_{34}^\prime$};
    \node[scale=1.0] at (.4, .4) {$\eb_{14}^\prime$};
    \node[scale=1.0] at (-.4, -.4) {$\eb_{41}^\prime$};

    \node[main node] (5) at (7, 2)   {$\vb_1$};
    \node[main node] (6) at (11, 2)   {$\vb_2$};
    \node[main node] (7) at (7, -2) {$\vb_3$};
    \node[color4, main node] (8) at (11, -2)  {$\vb_4^\prime$};

    \path[every node/.style={font=\sffamily\small}]
    (5) edge[black] (6)
    (5) edge[black] (7)
    (5) edge[black, bend left] (8)
    (6) edge[black] (8)
    (7) edge[black] (8)
    (8) edge[black, bend left] (5);

    \node[scale=1.0] at (9, 2.5) {$\eb_{12}^\prime$};
    \node[scale=1.0] at (6.5, 0) {$\eb_{13}^\prime$};
    \node[scale=1.0] at (11.6, 0) {$\eb_{24}^\prime$};
    \node[scale=1.0] at (9, -2.5) {$\eb_{34}^\prime$};
    \node[scale=1.0] at (9.4, .4) {$\eb_{14}^\prime$};
    \node[scale=1.0] at (8.6, -.4) {$\eb_{41}^\prime$};

    \draw[-{Triangle[width=18pt, length=8pt]}, line width=18pt] (3.5,0) -- (5.7,0);
    \draw[-{Triangle[width=18pt, length=8pt]}, line width=18pt] (3.5-9,0) -- (5.7-9,0);

    \node[white, scale=1.0] at (4.5-9, 0) {$\rho_{e \rightarrow v}$};
    \node[white, scale=1.0] at (4.5, 0) {$\phi_v$};


    \node[scale=1.0] at (-4.5, 5) {$\displaystyle \rho_{e \rightarrow v}\left( \textcolor{color2}{\eb_{34}^\prime}, \textcolor{color1}{\eb_{24}^\prime}, \textcolor{color3}{\eb_{14}^\prime} \right) = \textcolor{color5}{\mathbf{\overline{e}}_4}$};
    \node[scale=.8] at (4.5, 5) {$\displaystyle \phi_{v}\left( \begin{bmatrix} \textcolor{teal}{\vb_4} \\ \textcolor{color5}{\mathbf{\overline{e}}_4} \\ \gb \end{bmatrix}\right) = \textcolor{color4}{\vb_4^\prime}$};

  \end{tikzpicture}
\end{figure}
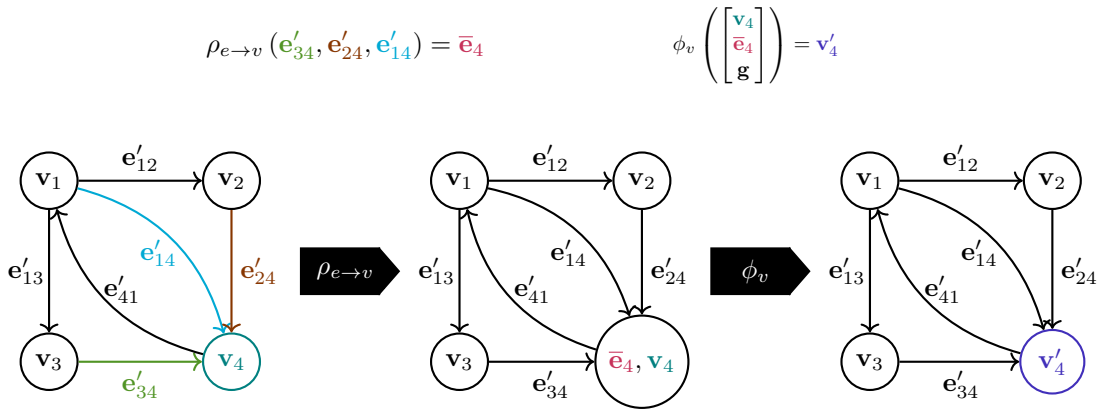

With the edges transformed, the next {\bf for} loop in Algorithm~\ref{algo:GNBlock} updates the
vertex attributes. Each iteration contains two steps, the first being the aggregation of all edges terminating
at the current vertex using $\rho_{e\rightarrow v}$, and the second updating the vertex's attribute.
This is depicted in Figure~\ref{fig:vertex-update} for $v_4$. Here, the edge neighborhood
contains edges that terminate at $v_4$. The attributes $\eb_{24}', \eb_{34}'$ and $\eb_{14}'$ 
for these three edges are aggregated using $\rho_{e\rightarrow v}$ yielding $\bar{\eb}_4$.
The aggregated edge attributes, with the vertex attribute $\vb_4$, and the global attribute ${\bf g}$
are input to the update function $\phi_v$. 
While only one update is shown, this procedure is repeated (potentially in parallel) until all 
vertex attributes are updated. 

The final phase of the GNN layer is to transform the global attributes of the graph. This is done by
first aggregating all the transformed edge attributes using $\rho_{e \rightarrow g}$ 
and all the transformed vertex attributes using $\rho_{v\rightarrow g}$  .
Then the new global attributes $\gb'$ are computed using 
$\phi_g$. The transformed edge, vertex, and graph attributes are returned in the
final line of Algorithm~\ref{algo:GNBlock}.

A notable quality of the message passing layer is that if the graph attributes are
ignored ($\gb$ doesn't change), then the action of the layer only considers the one-ring
(or distance-one)
neighborhood of a vertex for attribute update. So the action of two consecutive GNN layers will
update attributes based on vertices two edges away, or the two-ring neighborhood. This breadth-first
approach has implications on the use of GNNs for processing subgraph information independently,
leaving intriguing possibilities for parallelism. Additionally, each pass of the layer results in
an exchange of information over attributes within a single vertex one-ring neighborhood. Consequently,
multiple applications of the layer communicate to
ever broader vertex neighborhoods. 

\paragraph{Viral Interlude}
The update for the viral outbreak graph described conceptually above can be made (more) formal by defining appropriate aggregation and update functions. The edge update function $\phi_e$ must provide a measurement of an interaction contributing to the cure and disease state of an individual. Further, if a graph function is to be applied multiple times then the time of interaction must also be represented in the results of this function:
\begin{equation}
{\bf e}_{ij}' \leftarrow \phi_e({\bf e}_{ij},{\bf v}_i,{\bf v}_j,{\bf g}) = [{\bf e}_{ij}, {\bf i}_{ij}]
\end{equation}
where ${\bf e}_{ij}$ is the original attribute, the time of interaction, and ${\bf i}_{ij}$ is a quantification of the impact of the interaction (note that in the context of neural networks the quantity ${\bf i}_{ij}$ is not understood in a precise way, rather the GNN will \emph{learn} to interpret it). For this layer $n_e' = n_e +1$. While this quantifies the impact of a single interaction, the end state of an individual is based on all their interactions specified by connected graph edges. This is the intent of aggregation function, $\rho_{e\rightarrow v}$.  For example, the sum of the interaction risks might be a useful metric, as would the maximum and/or minimum risk, depending on the dynamics of the cure and disease. More formally, a vector-valued function that quantifies the external risk to an individual $j$ over all interactions is
\begin{equation}
\bar{\bf e}_j = \rho_{e\rightarrow v}({\bf e}_{*j}') = \left[\sum_{i} {\bf e}_{ij}'~~, \max_i {\bf e}_{ij}',~~ \min_i {\bf e}_{ij}'\right].
\end{equation}
With the external risks to an individual quantified, the vertex attributes specifying the cure/disease probability for the individual must be updated using $\phi_v$.  

\subsection{Training}\label{sec:GNN-training}
Typically in GNN applications, MLP neural networks are used as the update functions. Each of
these MLPs is parameterized by a set of weights and biases. Therefore, the parameters of
a GNN is the union of the weights and biases of all the update functions. 
We denote the collected set of parameters for the neural network as $\Theta$. The activity of tuning these parameters
to learn to perform a set task is referred to as training. The are different types of tasks that can be learned, but
the most common is referred to as supervised learning. In this case, a desired mapping is specified
empirically using a dataset composed of input-output pairs. For exposition, in this section we define the data
set as $(X_k,Y_k)$ for $k=1\ldots n_d$, where $X_k$ is the domain element mapping to $Y_k$, the range element.
For GNN's, the domain and range values are defined as vertex, edge and/or global attributes. 
The parameters of a neural network are said to be trained when a \emph{loss} function is minimized. An idealized
supervised learning problem and loss function is
\begin{equation}\label{eqn:ideal-loss}
\Theta^* = \argmin_\Theta  \left(\mathcal{L}(\Theta) := \frac{1}{n_d}\sum_{k=1}^{n_d} {l}_{\Theta}(\mathcal{N}(X_k;\Theta),Y_k)\right)
\end{equation}
where the function $l_\Theta$ measures the difference between the GNN prediction for the attributes $X_k$ and the empirical target $Y_k$.
Choices for $l_\Theta$ include norms for regression, or cross entropy for classification. 


With the loss function selected, 
an algorithmic approach to find the optimal parameters $\Theta^*$ is required. 
Initially the weights and biases are selected at random (see~\cite{glorot2010understanding,he2015delving,goodfellow2016deep,cyr2020robust} for example approaches and considerations).
Then an iterative method is used to incrementally reduce the loss function and improve the prediction by adjusting
the GNN parameters.
The most common iterative
optimization methods for training neural networks are gradient descent algorithms. At each iteration, 
the algorithm calculates the loss for the current $\Theta$ parameters by forward propagation through
a GNN and evaluating a loss on the predicted versus observed output values. In a second step of the iteration, the gradient of the loss with respect to $\Theta$ (e.g. $\nabla_\Theta \mathcal{L}$) is computed using the celebrated \emph{backpropagation} algorithm~\cite{goodfellow2016deep}. 
Simplistically, the negative of the gradient is used to update the parameters and the iteration is repeated.
An important metric to measure the rate of the reduction of the loss is the \emph{epoch}. One epoch corresponds to each
entry in the dataset having been used in a gradient descent step once. Values of tens to thousands of epochs are not uncommon. 
While gradient descent~\cite{nocedal1999numerical} is a common algorithm in training,  
more sophisticated algorithms can lead to more robust results, for example
Stochastic Gradient Descent~\cite{robbins1951stochastic},
RMSProp~\cite{tieleman2012lecture}, and Adam~\cite{Adam14}. For a review of optimization methods
for machine learning see~\cite{bottou2018optimization,goodfellow2016deep}.

The ability to calculate the gradient for a large number of parameters is automated
by software libraries which perform automatic differentiation. Common libraries for
machine learning such as PyTorch~\cite{PyTorch2015}, TensorFlow~\cite{tensorflow2015-whitepaper},
Jax~\cite{jax2018github}, and others all include this capability, allowing gradient descent
algorithms to be applied without requiring user-specified gradient functions.

\paragraph{Viral Interlude} Applied to our viral example, the goal of training is to determine parameters that define
the transmission rates, and effectiveness of the cure versus the spread of the disease. In this supervised learning example,
the dataset input is defined as the initial distribution of the virus
and cure over the vertices of the graph, while the output data
will be the distribution of the virus and cure after multiple cycles. The dataset would be comprised of graphs
defined by multiple communities already 
inflicted by the virus (thus an observation can be obtained). 
The model edge update function is parameterized
by a neural network that defines the transfer of cure/disease through an individual interaction. The vertex update is also
a neural network model that defines the uptake of the cure/disease by a single individual through multiple interactions. Prior
to training, the neural network contains only the context of the interactions, not the specifics of the cure/viral diffusion. These
specifics would be learned by defining a regression-based loss on the vertex attributes associated with the probability of an
individual having the virus/cure. Once trained, the GNN model can be applied to communities (graphs) that have an initial
viral load to predict a final state distribution.

\section{Linear Operations as Graph Network Layers}\label{sec:lin_gnn}
This section demonstrates how linear algebra computations can be
carried out utilizing the language and structure of graph neural networks.
These familiar algorithms are selected to make the
graph neural network framework described above concrete.
We note that the representation of these
operations as graph neural networks are not unique, and multiple
representations may exist for a given calculation. As such, this
is not meant to be an exhaustive list, but instead a general guide
to setting up and representing different linear algebra operations as
graph neural networks.

We begin by briefly illustrating the natural relationship between sparse matrices
and graphs. A weighted graph with $n$ vertices 
can represent a square $n \times n $ matrix $A$. Each 
nonzero $A_{ij}$ defines a weight for a directed edge that
emanates from vertex $j$ and terminates at vertex $i$.
Notice that this definition uses self edges to represent 
entries on the matrix diagonal. Alternatively, one can 
omit self edges and instead store $A_{ii}$ at the $i^{th}$ vertex. 
The inclusion or omission of self edges affects how linear algebra 
operations are represented and how data  propagates through
the network as different information is available during specific 
update/aggregation phases.  The algorithms in this section assume that
the graph neural network includes self edges, except for 
the matrix-vector example in Section~\ref{sec:algo_matvec}.
Additionally, our edge orientation choice
is natural for matrix-vector products
$ y = A x $ where the $i^{th}$ entry of 
$y$ is defined by
\begin{equation*}
  y_{i} = \sum_{j \in A_{i*}} A_{ij}x_{j} 
\end{equation*}
and $A_{i*}$ denotes the set of all nonzero entries in the $i^{th}$ matrix row.
Here, information from neighboring edges must be gathered at the $i^{th}$ vertex.
This is accomplished by an aggregation function as the edges in the
$i^{th}$ matrix row terminate at the $i^{th}$ vertex.
More generally, our edge orientation choice emphasizes that information 
that directly influences vertex $i$ flows into vertex $i$.
Finally, we note that we only consider square matrices. However, 
non-square matrices can be represented as bipartite graphs where rows 
and columns each have their own distinct set of vertices. 
Bipartite graphs extensions of Algorithm~\ref{algo:GNBlock} are possible, but
are beyond the scope of this educational survey.  

We start with foundational linear algebra algorithms:
sparse matrix-vector product and a matrix-weighted norm. 
Next,
we describe three simple iterative methods: a weighted Jacobi linear solver,
a Chebyshev linear solver, and a power method eigensolver. We conclude
with some kernels used within an advanced algebraic multigrid
linear solver. 
All of these examples follow
Algorithm~\ref{algo:GNBlock} for the structure of the graph neural
network. The details of each operation specific component, such as the
update functions ($\phi_v$, $\phi_e$, $\phi_g$) and aggregation
functions ($\rho_{e\rightarrow v}$, $\rho_{e\rightarrow g}$,
$\rho_{v\rightarrow g}$), are summarized in
Tables~\ref{algo:spmv}-\ref{algo:soc_classic}.

\paragraph{\textbf{Code Availability}}
Code for each of the GNNs given in this chapter can be found in the repository located at \url{https://github.com/sandialabs/gnn-applied-linear-algebra/}. Code is provided both in MATLAB script and python via PyTorch and the PyTorch
Geometric package. Each example provides the code implementation of the layer as well as a small
demonstration which shows the output of the GNN matches that of the ``traditional'' method.


\newcommand{\GNNDataAndLayer}[9]{%
  \def\GNNDataAndLayerMore##1##2##3##4{%
            &  Fixed       &  Mutable         &  \multicolumn{2}{|c|}{Updates} &%
     \multicolumn{2}{|c|}{Aggregation} \\%
     \hline\hline%
            &    &    &  \multicolumn{4}{|c|}{#7} \\
     \cline{4-7}%
     Edge   & #1 & #2 & $\phi_e$ & #8  & $\rho_{e\rightarrow v}$  & ##2 \\%
     Vertex & #3 & #4 & $\phi_v$ & #9  & $\rho_{e\rightarrow g}$  & ##3 \\%
     Global & #5 & #6 & $\phi_g$ & ##1 & $\rho_{v\rightarrow g}$  & ##4 \\%
  }
  \GNNDataAndLayerMore
}

\newcommand{\GNNExtraLayer}[7]{%
  \cline{4-7}%
  & & & \multicolumn{4}{|c|}{#1} \\
  \cline{4-7}%
  & & & $\phi_e$ & #2  & $\rho_{e\rightarrow v}$ & #5 \\%
  & & & $\phi_v$ & #3  & $\rho_{e\rightarrow g}$ & #6 \\%
  & & & $\phi_g$ & #4  & $\rho_{v\rightarrow g}$ & #7 \\%
}


\newenvironment{gnntable}
{\begin{center}
\begin{tabular}{|lll||ll|ll|}%
  \hline%
  \multicolumn{3}{|c||}{Data} & \multicolumn{4}{|c|}{Functions} \\%
}
{\hline
  \end{tabular}
  \end{center}
}
\subsection{Sparse Matrix-Vector Product}\label{sec:algo_matvec}

Sparse matrix-vector products are fundamental building blocks of
many numerical linear algebra algorithms, such as Krylov methods.
Table~\ref{algo:spmv} illustrates a GNN that computes $y=Ax$.
\begin{table}[htb]
  \caption{Sparse Matrix-Vector Product as a Graph Network with Self-Edges}\label{algo:spmv}
  \begin{gnntable}
    \GNNDataAndLayer
    {$A_{ij}$}{$c_{ij}$}
    {$x_i$}{$y_i$ [Output]}
    {---}{---}
    {Layer 1}
    {$c_{ij} = A_{ij} x_j$}{$y_i=\overline{c}_i$}{---}
    {$\overline{c}_{i} = \sum_j c_{ij}$}{---}{---}
  \end{gnntable}
\end{table}
The left block of the table describes the input, intermediate and 
output data for each graph component (vertex, edge, global). The 
right column is subdivided into the update functions required by
Algorithm~\ref{algo:GNBlock} on the left, and the aggregation functions.
Further, if there are distinct layers in the network, the right column
lists the sets of update and aggregation functions. 

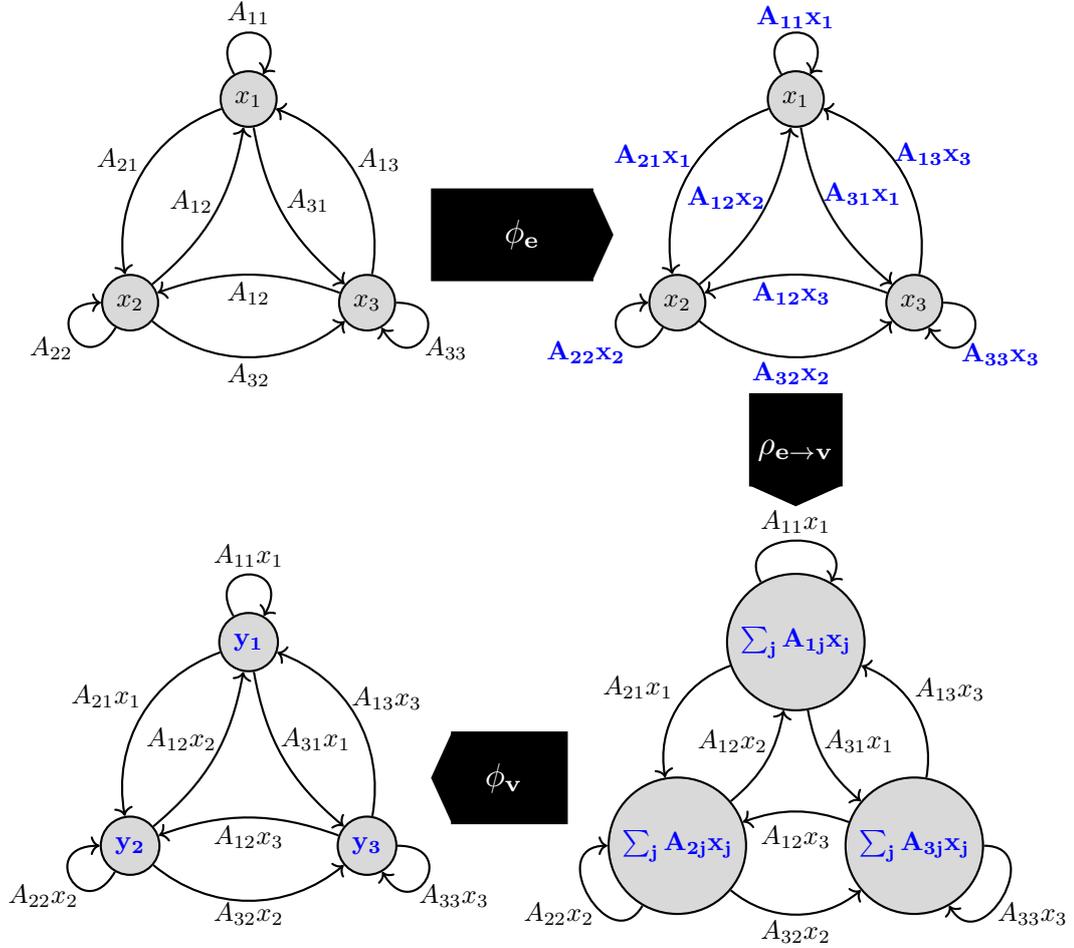
\begin{figure}
  \caption{A graphical depiction of a GNN layer described by Table~\ref{algo:spmv} that shows the state of the
  GNN and transitions defined by Algo.~\ref{algo:GNBlock}.}
    \label{fig:spmv_gnn_graphs}
\vskip .1in
  \centering
  \begin{tikzpicture}[->,auto,
    thick,main node/.style={circle,draw,font=\sffamily\normalsize,fill=black!15}, scale=.6]
    \node[main node] (1) at (90:3)  {$x_1$};
    \node[main node] (2) at (210:3) {$x_2$};
    \node[main node] (3) at (330:3) {$x_3$};

    \path[every node/.style={font=\sffamily\small}]
    (1) edge[black, loop, in=60, out=120, min distance=40] (1)
    (1) edge[black, bend right=40] (2)
    (1) edge[black, bend right=20] (3)
    (2) edge[black, bend right=20] (1)
    (2) edge[black, loop, in=180, out=240, min distance=40] (2)
    (2) edge[black, bend right=40] (3)
    (3) edge[black, bend right=40] (1)
    (3) edge[black, bend right=20] (2)
    (3) edge[black, loop, in=300, out=360, min distance=40] (3);

    \node[scale=1.0] at (90:4.9)  {$A_{11}$};
    \node[scale=1.0] at (150:3.3) {$A_{21}$};
    \node[scale=1.0] at (30:1.5)  {$A_{31}$};
    \node[scale=1.0] at (150:1.45) {$A_{12}$};
    \node[scale=1.0] at (210:5)   {$A_{22}$};
    \node[scale=1.0] at (270:3.1)   {$A_{32}$};
    \node[scale=1.0] at (30:3.3)  {$A_{13}$};
    \node[scale=1.0] at (270:1.3) {$A_{12}$};
    \node[scale=1.0] at (330:5)   {$A_{33}$};

    \begin{scope}[shift={(12, 0)}]
      
      \node[main node] (1) at (90:3)  {$x_1$};
      \node[main node] (2) at (210:3) {$x_2$};
      \node[main node] (3) at (330:3) {$x_3$};

      \path[every node/.style={font=\sffamily\small}]
      (1) edge[black, loop, in=60, out=120, min distance=40] (1)
      (1) edge[black, bend right=40] (2)
      (1) edge[black, bend right=20] (3)
      (2) edge[black, bend right=20] (1)
      (2) edge[black, loop, in=180, out=240, min distance=40] (2)
      (2) edge[black, bend right=40] (3)
      (3) edge[black, bend right=40] (1)
      (3) edge[black, bend right=20] (2)
      (3) edge[black, loop, in=300, out=360, min distance=40] (3);

      \node[scale=1.0] at (90:4.8)  {\textcolor{blue}{ $\mathbf{A_{11}x_1}$}};
      \node[scale=1.0] at (150:3.5) {\textcolor{blue}{ $\mathbf{A_{21}x_1}$ }};
      \node[scale=1.0] at (30:1.8)  {\textcolor{blue}{ $\mathbf{A_{31}x_1}$ }};
      \node[scale=1.0] at (150:1.65) {\textcolor{blue}{ $\mathbf{A_{12}x_2}$ }};
      \node[scale=1.0] at (210:5.2) {\textcolor{blue}{ $\mathbf{A_{22}x_2}$ }};
      \node[scale=1.0] at (270:3.05)   {\textcolor{blue}{ $\mathbf{A_{32}x_2}$ }};
      \node[scale=1.0] at (30:3.6)  {\textcolor{blue}{ $\mathbf{A_{13}x_3}$ }};
      \node[scale=1.0] at (270:1.3) {\textcolor{blue}{ $\mathbf{A_{12}x_3}$ }};
      \node[scale=1.0] at (330:5.3) {\textcolor{blue}{ $\mathbf{A_{33}x_3}$ }};

    \end{scope}
    \begin{scope}[shift={(12, -12)}]
      
      \node[main node] (1) at (90:3)  {\textcolor{blue}{ $\mathbf{\sum_{j} A_{1j}x_{j}}$} };
      \node[main node] (2) at (210:3) {\textcolor{blue}{ $\mathbf{\sum_{j} A_{2j}x_{j}}$} };
      \node[main node] (3) at (330:3) {\textcolor{blue}{ $\mathbf{\sum_{j} A_{3j}x_{j}}$} };

      \path[every node/.style={font=\sffamily\small}]
      (1) edge[black, loop, in=60, out=120, min distance=40] (1)
      (1) edge[black, bend right=40] (2)
      (1) edge[black, bend right=20] (3)
      (2) edge[black, bend right=20] (1)
      (2) edge[black, loop, in=180, out=240, min distance=40] (2)
      (2) edge[black, bend right=40] (3)
      (3) edge[black, bend right=40] (1)
      (3) edge[black, bend right=20] (2)
      (3) edge[black, loop, in=300, out=360, min distance=40] (3);

      \node[scale=1.0] at (90:5.65)  {$A_{11}x_1$};
      \node[scale=1.0] at (150:4.0) {$A_{21}x_1$};
      \node[scale=1.0] at (30:1.6)  {$A_{31}x_1$};
      \node[scale=1.0] at (150:1.6) {$A_{12}x_2$};
      \node[scale=1.0] at (210:6.0) {$A_{22}x_2$};
      \node[scale=1.0] at (270:3.4) {$A_{32}x_2$};
      \node[scale=1.0] at (30:3.9)  {$A_{13}x_3$};
      \node[scale=1.0] at (270:1.3) {$A_{12}x_3$};
      \node[scale=1.0] at (330:6.0) {$A_{33}x_3$};

    \end{scope}

    \begin{scope}[shift={(0, -12)}]
      
      \node[main node] (1) at (90:3)  {\textcolor{blue}{$\mathbf{y_1}$} };
      \node[main node] (2) at (210:3) {\textcolor{blue}{$\mathbf{y_2}$} };
      \node[main node] (3) at (330:3) {\textcolor{blue}{$\mathbf{y_3}$} };

      \path[every node/.style={font=\sffamily\small}]
      (1) edge[black, loop, in=60, out=120, min distance=40] (1)
      (1) edge[black, bend right=40] (2)
      (1) edge[black, bend right=20] (3)
      (2) edge[black, bend right=20] (1)
      (2) edge[black, loop, in=180, out=240, min distance=40] (2)
      (2) edge[black, bend right=40] (3)
      (3) edge[black, bend right=40] (1)
      (3) edge[black, bend right=20] (2)
      (3) edge[black, loop, in=300, out=360, min distance=40] (3);

      \node[scale=1.0] at (90:4.9)  {$A_{11}x_1$};
      \node[scale=1.0] at (150:3.6) {$A_{21}x_1$};
      \node[scale=1.0] at (30:1.7)  {$A_{31}x_1$};
      \node[scale=1.0] at (150:1.7) {$A_{12}x_2$};
      \node[scale=1.0] at (210:5.2) {$A_{22}x_2$};
      \node[scale=1.0] at (270:3)   {$A_{32}x_2$};
      \node[scale=1.0] at (30:3.5)  {$A_{13}x_3$};
      \node[scale=1.0] at (270:1.3) {$A_{12}x_3$};
      \node[scale=1.0] at (330:5.2) {$A_{33}x_3$};

    \end{scope}

    \draw[-{Triangle[width=35pt, length=8pt]}, line width=35pt] (4,0) -- (8,0);
    \draw[-{Triangle[width=35pt, length=8pt]}, line width=35pt] (12,-3.5) -- (12,-6);
    \draw[-{Triangle[width=35pt, length=8pt]}, line width=35pt] (7,-12) -- (4,-12);

    \node[white, scale=1.25] at (6, 0.0) {$\mathbf{\phi_{e}}$};
    \node[white, scale=1.25] at (12, -4.75) {$\mathbf{\rho_{e \rightarrow v}}$};
    \node[white, scale=1.25] at (5.6, -12) {$\mathbf{\phi_{v}}$};

  \end{tikzpicture}
\end{figure}

For the sparse matrix vector product the nonzero entries $A_{ij}$ naturally
correspond to edges, so we assign them to edges as fixed (i.e., unmodified by the GNN) data objects.
The vector entries $x_i$ and $y_i$ are assigned to vertices, with $x_i$ fixed and $y_i$ mutable (i.e., modified by the GNN) and
initialized to zero.  The upper left graph in Figure~\ref{fig:spmv_gnn_graphs} illustrates the starting condition
of the GNN with data placed according to the description in Table~\ref{algo:spmv}.

We now follow the update and aggregation
functions of Algorithm~\ref{algo:GNBlock} in order. 
Each step is a transition between graph states in Figure~\ref{fig:spmv_gnn_graphs}.  
The edge feature update $\phi_e$ takes fixed edge feature $A_{ij}$
and fixed vertex feature $x_j$ and multiplies them, storing the result on
each edge.
Summation is employed for edge aggregation, denoted $\phi_e$,  with 
the result stored on the vertex in the vertex update phase, $\phi_v$. The lower left
graph in Figure~\ref{fig:spmv_gnn_graphs} shows the final state of the GNN, with the
entries $y_i$ containing the matrix vector product located at graph vertices.
In Table ~\ref{algo:spmv} and all the tables that follow, it is understood that operations
such as $A_{ij} x_j$ and $\sum_j c_{ij}$ only operate on nonzero entries of 
$A$ and $c$ respectively. That is, the edges must be included in the graph.
The table can be modified if self-edges are not stored. To that end, Table~\ref{algo:spmv2} shows a slightly more
complicated GNN.  The difference is that
$A_{ii}$  is stored on vertices and so the $A_{ii} x_i$ term is added in 
the vertex update function $\phi_v$.

\begin{table}[htb]
  \caption{Sparse Matrix-Vector Product as a Graph Network without Self-Edges}\label{algo:spmv2}
  \begin{gnntable}
    \GNNDataAndLayer
    {$A_{ij}$}{$c_{ij}$}
    {$x_i,A_{ii}$}{$y_i$ [Output]}
    {---}{---}
    {Layer 1}
    {$c_{ij} = A_{ij} x_j$} {$y_i=\overline{c}_i + A_{ii} x_i$} {---}
    {$\overline{c}_i = \sum_j c_{ij}$}{---}{---}
  \end{gnntable}
\end{table}

\subsection{Matrix-Weighted Norm}\label{sec:algo_wnorm}

Similar to the matrix-vector product, a GNN for computing a matrix-weighted vector norm, $\norm{x}_W
= \sqrt{x^T W x}$, is shown in Table~\ref{algo:mat_weighted_norm}. 
\begin{table}[htb!]
  \caption{Matrix-Weighted Norm}\label{algo:mat_weighted_norm}
  \begin{gnntable}
    \GNNDataAndLayer
    {$W_{ij}$}{$c_{ij}$}
    {$x_i$}{$y_i$}
    {---}{$n$ [Output]}
    {Layer 1}
    {$c_{ij} = W_{ij} x_j$} {$y_i= x_i\overline{c}_i$} {$n = \sqrt{\overline{y}}$}
    {$\overline{c}_i = \sum_j c_{ij}$}{---}{$\overline{y} = \sum_i y_i$}
  \end{gnntable}
\end{table}
The matrix-vector product from
Table~\ref{algo:spmv} is modified in the vertex update and adds operations for
vertex-to-global aggregation and global update. After the edge update and the
edge-to-vertex aggregation, $\overline{c}$ contains the vector $Wx$. Therefore,
in the vertex update, we multiply by $x_{i}$ to obtain the vector $y_{i}$. In
the vertex-to-global aggregation, $y_{i}$ is summed to obtain $x^{T} W x$.
Finally, in the global update, the square root is taken to yield the vector
norm.
\subsection{Weighted Jacobi Iteration}\label{sec:algo_jacobi}
Weighted Jacobi iteration is a 
simple method to solve a linear system $A x = b$ for $x$.  
One iteration of weighted Jacobi is written as an update formula
\begin{equation}\label{eq:jacobi}
x^{k+1} = x^{k} + \omega D^{-1}(b-A x^{k}),
\end{equation}
where $k$ is the iteration index, $D$ is the matrix diagonal, $b$ is the right-hand side,
$x^{k}$ is the solution at the $k^{th}$ Jacobi iteration, and
$\omega$ is the weight parameter. For the GNN representation of the weighted Jacobi
method, the edge data includes the matrix entries $A_{ij}$.  The initial vertex data is comprised
of the matrix diagonal,
$A_{ii}$, and the right-hand side vector, $b_i$. The weight parameter, $\omega$ is included as
fixed global data. Note that while the matrix diagonal, $A_{ii}$, is stored
as vertex data, this GNN layer uses self-edges as well, so the matrix diagonal is
also accessible as regular edge data. This is done to allow easy access to this
data for both the matrix-vector product and the Jacobi update. For the functions, all of the
components of the sparse matrix-vector product are retained, except for the
vector update, $\phi_v$. That is replaced to use the aggregated edge data inside a
Jacobi-style update, \eqref{eq:jacobi} to yield the GNN layer shown in
Table~\ref{algo:jacobi}.
\begin{table}[h!]
\caption{Weighted Jacobi Iteration as a Graph Network}\label{algo:jacobi}
\begin{gnntable}
  \GNNDataAndLayer
  {$A_{ij}$}{$c_{ij}$}
  {$A_{ii}, b_i$}{$x_i$ [Output]}
  {$\omega$}{---}
  {Layer 1}
  {$c_{ij} = A_{ij} x_j$} {$x_i=x_i + \omega A_{ii}^{-1} (b_i-\overline{c}_i)$}{---}
  {$\overline{c}_i = \sum_j c_{ij}$}{---}{---}
\end{gnntable}
\end{table}
\subsection{Chebyshev Iterative Solver}\label{sec:algo_cheby}
The Chebyshev method is another simple iterative scheme.
It can be viewed as a type of multi-step Jacobi algorithm where different
weights are used for each step and these weights are optimal in some sense given
knowledge about the largest and smallest eigenvalue ($\lambda_{max}$ and $\lambda_{min}$)
of the symmetric positive definite matrix $A$.
An example algorithm is given in
Algorithm~\ref{algo:pseudo_cheby}, following 
the algorithm presentation in \cite{saadBook}. 

\begin{algorithm}
\caption{Chebyshev Solver}\label{algo:pseudo_cheby}
\begin{algorithmic}[0]
\Require Matrix $A$, right-hand side $b$, initial guess $x_0$, $\#$ of iterations $N$, $\lambda_{max}$ \& $\lambda_{min}$
\Ensure Solution $x$

\\\hrulefill

  \Function{$(x)$ = Cheby}{$A,b,x_0,\lambda_{max},\lambda_{min}$,$N$}

  \State $x \gets x_0$ ;\hskip 1.01in $r \gets b - A x$
  \State $\theta \gets ( \lambda_{max} + \lambda_{min} )/2$ ;~~ $\delta \gets ( \lambda_{max} - \lambda_{min} )/2$ ;~~$\sigma \gets \theta/\delta$ 
  \State $\rho \gets 1/\sigma $ ; \hskip .88in $d \gets (1/\theta)~ r$

  \For{$i \in \{1\ldots N\}$}
    \State $x \gets x + d$
    \State $r \gets r - A d $
    \State $\rho_{\text{prior}} \gets \rho$
    \State $\rho \gets 1/(2\sigma - \rho)$
    \State $d    \gets \rho~ \rho_{\text{prior}}~ d ~+~ (2 \rho/\delta)~r $
  \EndFor
  \State \Return $x$
  \EndFunction
\vspace*{-.02in} 
\end{algorithmic}
\end{algorithm}

\begin{table}[h!]
   \caption{Chebyshev Solver Iteration as a Graph Network}\label{algo:cheby} 
   \begin{gnntable}
      \GNNDataAndLayer
      {$A_{ij}$}{$c_{ij}$}
      {}{$r_{i}, d_{i}, x_{i}$ [Output] }
      {$\delta$, $\sigma$}{$\rho, \rho_{\text{prior}} $}
      {Layer 1}
      {---}{$x_{i} = x_{i} + d_{i}$}{---}
      {---}{---}{---}
      \GNNExtraLayer{Layer 2}
      {$c_{ij} = A_{ij} d_{j}$}{$r_{i} = r_{i} - \overline{c}_{i}$ }{$\rho_{\text{prior}}  = \rho ; \rho = 1 / (2\sigma - \rho)$}
      {$\overline{c}_{i} \sum_{j} c_{ij}$ }{---}{---} 
      \GNNExtraLayer{Layer 3}
      {}{$d_{i} = \rho~ \rho_{\text{prior}}~ d_{i}  ~+~ (2 \rho/\delta)~r_{i}$}{}
      {}{}{}
   \end{gnntable}
\end{table}

The GNN describing one iteration of the Chebyshev method is outlined in Table~\ref{algo:cheby},
where all of the initial data values are assumed to have the specified values assigned before the
{\bf for} loop in Algorithm~\ref{algo:pseudo_cheby}. For the Chebyshev GNN, most of the computation
is nodal and the GNN primarily serves to propagate residual information. When multiple values are
updated on a single line of an update function, the values update from left to right, so parameters
are updated in the correct order. In each iteration of the Chebyshev algorithm the GNN uses Layer 1
to update the current iterate. Then Layer 2 updates the residual information. Finally, Layer 3
updates the direction vector based on the the updated residual. If $N$ iterations of the Chebyshev
method are desired, run the GNN in Table~\ref{algo:cheby} $N$ times consecutively. 
\subsection{Power Method}\label{sec:algo_power_method}
Since the power method consists for symmetric matrices consists of matrix-vector multiplications and
norms, we can represent the power method as a GNN as well. The details are given
in Tables~\ref{algo:power_method_iterate} and \ref{algo:power_method_rayleigh}.
We split the algorithm into two networks. The network in
Table~\ref{algo:power_method_iterate} is the iterative portion which updates the
\begin{table}[h!]
	\caption{Power Method as a Graph Network - Iterative Layers}\label{algo:power_method_iterate}
  \begin{gnntable}
    \GNNDataAndLayer
    {$A_{ij}$}{$c_{ij}$}
    {---}{$b_i, y_i$}
    {---}{$h_{1} , n_A, \lambda_{\text{max}} [Output]$}
    {Layer 1 ($b = Ab$)}
    {$c_{ij} = A_{ij}b_j$} {$b_i = \overline{c}_i$} {---}
    {$\overline{c}_i = \sum_j c_{ij}$}{---}{---}
    \GNNExtraLayer{Layer 2 ($h_1 = \norm{b}$)}
    {---}{$y_i = b_i^2$}{$h_{1}  = \sqrt{\overline{y}}$}
    {---}{---}{$\overline{y} = \sum_i y_i$}
    \GNNExtraLayer{Layer 3 (Re-Normalize)}
    {---}{$b_i = \frac{b_i}{h_{1} }$}{---}
    {---}{---}{---}
  \end{gnntable}
\end{table}
$b$ vector to an improved approximation of the eigenvector associated with the
dominant eigenvalue. The second network computes the Rayleigh quotient to give
the eigenvalue estimate. Thus, for more than a single iteration of the power
\begin{table}[h!]
	\caption{Power Method as a Graph Network - Rayleigh Quotient}\label{algo:power_method_rayleigh}
  \begin{gnntable}
    \GNNDataAndLayer
    {$A_{ij}$}{$c_{ij}$}
    {---}{$b_i, y_i$}
    {---}{$h_{1} , h_{2} , \lambda_{\text{max}} [Output]$}
    {Layer 1 ($h_{2}  = b^T A b$)}
    {$c_{ij} = A_{ij} b_j$}{$y_i= b_i \overline{c}$}{$h_{2}  = \overline{y}$}
    {$\overline{c}_i = \sum_j c_{ij}$}{---}{$\overline{y} = \sum_i y_i$}
    \GNNExtraLayer{Layer 2 (Rayleigh Quotient)}
    {---}{$y_i = b_i^2$}{$\lambda_{\text{max}} = \frac{h_{2}}{\overline{y}}$}
    {---}{---}{$\overline{y} = \sum_i y_i$}
  \end{gnntable}
\end{table}
method, the network in Table~\ref{algo:power_method_iterate} is appplied successively
for the desired number of iterations, and then the Rayleigh quotient is computed from
the network in Table~\ref{algo:power_method_rayleigh}. 

\subsection{Advanced Algebraic Mutligrid Solver}\label{sec:algo_AMG}
We conclude this section with two algebraic multigrid kernels: one for interpolation and another associated
with coarsening a matrix graph. Algebraic multigrid (AMG) is an advanced linear solver.
In the description that follows, we consider solving the matrix system $A_1 x_1 = b_1$ for the 
unknown vector $x_1$.  While Jacobi and Chebyshev
may require a very large number of iterations to reach an acceptable solution, AMG 
can converge rapidly on many important matrix systems.  For example, AMG is known to converge efficiently on matrices
arising from discrete approximations of elliptic partial differential equations (PDEs) such as the heat equation \cite{Briggs2000}. 
In particular, the number of required iterations can be independent of the dimension of the matrix
$A_1$. Thus, increasingly larger matrix systems do not necessarily require increasingly more iterations. As matrix
systems 
can include 
over a billion unknowns, this property is extremely 
advantageous.  AMG methods are primarily used to solve discrete versions of PDEs (including many non-elliptic
PDEs), though they have been applied to other systems such as those arising from 
certain network problems.  While a full understanding of AMG is beyond our scope (see \cite{Briggs2000} for a more thorough explanation of the method), we present a few AMG ideas using a simple example.

Consider a metal plate with a hole where the plate is clamped at two locations. One clamp is
fixed at $40^\degree{\hskip -.01in}C$ while the other clamp is fixed at $10^\degree{\hskip -.01in}C$.  Suppose a discrete representation of Poisson's equation is 
used to model the steady-state heat distribution within the plate where the associated matrix $A_1$ only has nonzeros on 
the diagonal and for each edge of a mesh.  The vector $b_1$ contains only a few nonzeros 
to represent the fixed temperature of the plate where it is held by the clamps.  The other boundary conditions indicate that the 
heat gradient normal to the plate boundaries is zero. When starting with an initial guess of $x_1^{(0)} = 0$, the first Jacobi 
iteration produces $x_1^{(1)} = \omega D_1^{-1} A_1 b_1$. Due to the sparsity pattern of $A_1$, $x_1^{(1)}$ must only have a few nonzeros,
all within a graph distance of one from one of the two clamp locations.  In general, the $k^{th}$ Jacobi iteration extends $x_1^{(k-1)}$'s 
nonzero regions along edges adjacent to these regions. Thus, one can see that many iterations are required for the influence of the 
two clamps to propagate throughout the entire plate, and so it should not be surprising that many further iterations are required to 
reach a converged solution.  As the mesh is refined, more iterations are needed to propagate information from the clamps to the rest
of the plate.  The situation is similar for $N$-step Chebyshev solvers, which propagate information $N$ times faster than Jacobi but are
also $N$ times as expensive.

\begin{figure}[ht]
  \centering
  \caption{Left: sample mesh hierarchy (clamps shown on ${\cal G}_1$). Right: data movement within an AMG iteration (Jacobi iterations occur at tan circles and Gaussian elmination occurs at dark brown circle).}
   \includegraphics[trim = 0.1in 0.0in 4.2in 0.0in, clip, height = 5.8cm,width=6.3cm]{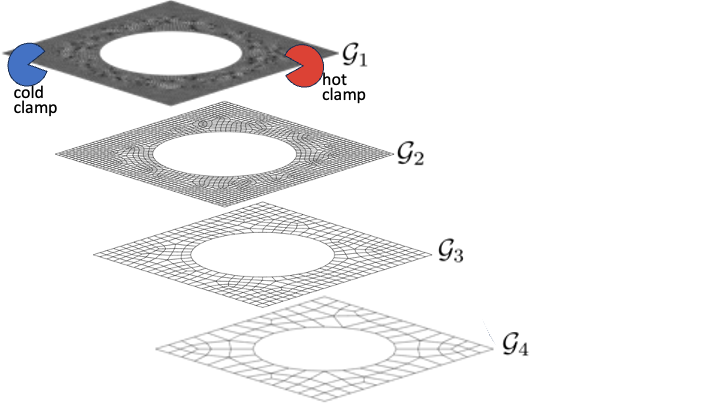}
   \includegraphics[trim = 6.6in 0.0in 0.1in 0.0in, clip, height = 5.8cm,width=6.3cm]{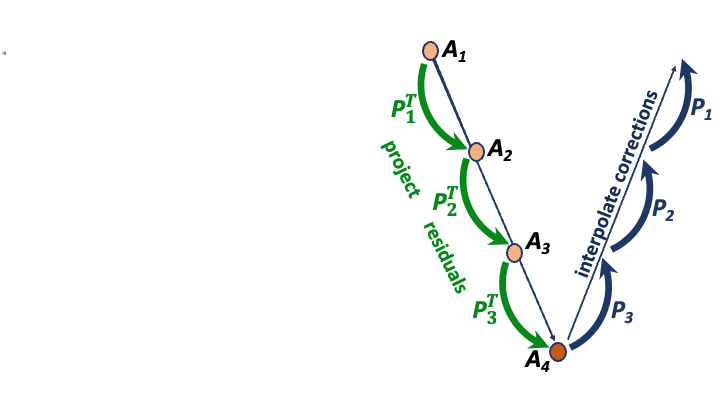}
  \label{fig:AMG Hierarchy}
\end{figure}

While Jacobi/Chebshev propagate information slowly across the mesh, it can be shown that these methods often efficiently  
produce approximate solutions where the error is smooth. The multigrid algorithm leverages this smoothing property in conjunction 
with a hierarchy of meshes ${\cal G}_\ell$ such as those depicted in Figure~\ref{fig:AMG Hierarchy}.  
One simple example of a multigrid iteration
employs $k$ Jacobi iterations on the $A_1 x_1 = b_1$ system to generate an approximate solution $x_1^{(k)}$. Since the error is now smooth,
it can be accurately approximated on a coarser mesh. Thus, the multigrid idea is to improve the solution using a 
coarser mesh. To do this, a residual is computed and then projected on to the ${\cal G}_2$ mesh within the hiearchy and used as the right 
hand side to a second linear system $A_2 x_2 = b_2$. Here, the matrix $A_2$ is a coarser discrete version of the Poisson equation 
constructed using ${\cal G}_2$. The $A_2$ matrix is a less accurate approximation to Poisson's equation,
but it is only used to correct the approximation obtained on the finer mesh. One can apply $k$ Jacobi iterations 
to this $A_2$ system to produce $x_2^{(k)}$. Notice that these iterations are less expensive than applying Jacobi to the $A_1$ system. 
We can repeat this same process (
apply Jacobi to $A_{\ell+1} x_{\ell+1} = b_{\ell+1}$ where $b_{\ell+1}$ is a projected version of the residual for the $A_\ell$ system) 
for all meshes in the hierarchy. In each case, the approximate solution on the $\ell^{th}$ mesh is a correction 
to the solution on the previous mesh system.  To complete the AMG iteration, we must add together the individual corrections.
That is, the approximate solution at the end of a single AMG iteration is given recursively by $\tilde{x}_\ell = x_\ell^{(k)} + P_{\ell} \tilde{x}_{\ell+1}$
where $P_{\ell}$ is a rectangular matrix that interpolates (or prolongates) approximate solutions associated with the $\ell{\hskip -.025in}+{\hskip -.025in}1$ mesh to
the $\ell$ mesh. Here $ 1 \le \ell \le N_{levels}-1$ and on the coarsest mesh we take $\tilde{x}_{N_{levels}} =  x_{N_{levels}}^{(k)}$.
Often a Gaussian elimination solver is used on the coarsest mesh as its cost is negligible when the matrix $A_{N_{levels}}$ is sufficiently
small.  Notice that in this case, one AMG iteration propagates information from the clamps throughout the entire mesh. While 
there is some cost to the coarse level computations, this cost is small relative to the $A_1$ computations if the coarser meshes
are significantly smaller than the finest mesh.
While it is generally difficult for application developers to generate the coarse operators ${\cal G}_\ell, A_\ell, $ and $P_\ell$ 
($P_\ell^T$ is commonly used for residual projections), an AMG algorithm automates this entire process.
That is, all of the additional operators/meshes are generated within the AMG method using graph algorithms that coarsen the 
matrix graph, define {\it algebraic} interpolation, and define the coarse $A_\ell$'s.

\subsubsection{Strength of Connection}\label{sec:algo_soc_sa}
Automatically coarsening a matrix graph is difficult. In our example, we need a coarse graph that can approximate
errors that remain after $k$ Jacobi steps. Several different AMG algorithms are used in practice each with its
own coarsening procedure. Typically, a strength-of-connection algorithm is applied as a pre-cursor to coarsening.
The basic idea is that some graph edges (or off-diagonal matrix nonzeros) may contribute very little to the 
information flow. For example, suppose one corner region of our plate 
consists of a near insulating material (e.g., extending from the corner to the closest segment of the circle). 
Edges between the conducting and insulating regions should be classified as weak so that they can effectively be ignored
during the graph coarsening process.  Unfortunately, improper edge classification can ruin AMG's impressive convergence 
properties and remains an active research topic that could possibly be improved by new ML algorithms.  We present 
two common strength-of-connection algorithms that are based on evaluating the relative size of off-diagonal matrix entries. 
However, we note that these algorithms may not be appropriate for complex matrix systems.  
%
Unlike the previously discussed methods that output information at graph vertices, 
strength of
connection algorithms output data on graph edges. Specifically, a strength of
connection algorithm produces a matrix $S(A)$ which has the same sparsity
pattern as the input matrix $A$. 

For the (non-symmetric) smoothed aggregation
AMG strength of connection, this computation is,
\begin{equation}\label{eq:soc_sa}
S_{ij} = \frac{A_{ij}^2}{A_{ii} A_{jj}}.
\end{equation}
The GNN layer version of this algorithm, which is shown in
Table~\ref{algo:soc_sa}, relies on fixed matrix edge data, $A_{ij}$, as well as
fixed diagonal values, $A_{ii}$, stored on the vertices. Note that this GNN
allows for self edges, so the diagonal values are stored on both the vertices and the edges.
The only calculation takes place in the edge update function,
$\phi_e$, which simply implements \eqref{eq:soc_sa}.

\begin{table}
  \caption{Smoothed Aggregation Strength of Connection as a Graph Network}\label{algo:soc_sa}
  \begin{gnntable}
    \GNNDataAndLayer
    {$A_{ij}$}{$S_{ij}$ [Output]}
    {$A_{ii}$}{---}
    {---}{---}
    {Layer 1}
    {$S_{ij} = A_{ij}^2 / (A_{ii} A_{jj})$}{---}{---}
    {---}{---}{---}
  \end{gnntable}
\end{table}

Similar to the smoothed aggregation strength of connection, the classical
strength of connection algorithm provides output data on graph edges to
construct a strength of connection matrix $S(A)$ that has the same sparsity as
the input matrix $A$. For the classic strength of connection, this computation
is,
\begin{equation}\label{eq:soc_classic}
    S_{ij} = \frac{-A_{ij}}{\max_{k\neq i}\{-A_{ik}\}}.
\end{equation}
The GNN representation of this algorithm, which is shown in
Table~\ref{algo:soc_classic}, relies only on fixed matrix edge data, $A_{ij}$.
The maximum negative value on each row is first determined in the aggregation
phase, $\rho_{e\rightarrow v}$, and is then propagated to the vertices. Finally, the edge
update is used to finish the calculation in~\eqref{eq:soc_classic}. 

One final step is needed to finish the classic strength of connection, and that is the
dropping of weak connections. For a suitable $0<\tau\leq 1$, drop the nonzero values
where $S_{ij}<\tau$. The sparsity pattern of $S_{ij}$ then provides the classic
strength of connection. This could be done as a post-processing step, or, to
keep with the theme of utilizing ideas from graph neural networks, we can modify
the final edge update to utilize a non-linear activation function, such as the \textbf{step} function: 
\vskip .01in
\begin{equation*}
  \text{step}(x) = \begin{cases}
    1 & x > 0 \\
    0 & x \leq 0
  \end{cases}\
\end{equation*}
to zero out the weak connections based on $\tau$. This results in
\begin{equation}\label{eq:soc_classic_drop}
    \hat{S}_{ij} = \text{step}\left(\frac{-A_{ij}}{\max_{k\neq i}\{-A_{ik}\}} - \tau \right),
\end{equation}
where the strong connections are the non-zero entries of $\hat{S}$.

\begin{table}
  \caption{Classic Strength of Connection as a Graph Network}\label{algo:soc_classic}
  \begin{gnntable}
    \GNNDataAndLayer
    {$A_{ij}$}{$S_{ij}$ [Output]}
    {---}{$v_i$}
    {---}{---}
    {Layer 1}
    {---}{$v_i = \overline{c}_i$}{---}
    {$\overline{c}_i = \max_{j\not=i} \{-A_{ij}\}$}{---}{---}
    \GNNExtraLayer{Layer 2}
    {$S_{ij} = -A_{ij} / v_i$}{---}{---}
    {---}{---}{---}
  \end{gnntable}
\end{table}

\subsubsection{Direct Interpolation}\label{sec:direct_interp}
We present one of the simplest AMG interpolation schemes, though more advanced/robust algorithms
are generally preferred. 
The interpolation algorithm provides output on the graph edges, similar to the strength of connection algorithms.
Specifically, we seek the nonzeros entries of $P_\ell$.  In the description that follows, 
we omit the sub-script $\ell$ to simplify notation.  The algorithm relies on first partitioning the vertex set $\mathcal{V}$ into 
$\mathcal{F}$ and $\mathcal{C}$ such that $\mathcal{F}\cap \mathcal{C} = \emptyset$ and $\mathcal{F}\cup \mathcal{C} = \mathcal{V}$
(cf., ~\cite{Briggs2000}).  The F-vertices only exist on the fine grid while C-vertices exist on both fine and coarse
grids.  We assume that the rows/columns of A (and hence the vertices) are ordered so that all C-vertices are numbered
before F-vertices.  In addition to partitioning, we require a strength of connection
matrix, $S$, as discussed in Section~\ref{sec:algo_soc_sa}. 

The direct interpolation operator is derived from input matrix
$A$ such that if $i$ is an F-vertex, then
\begin{equation}\label{eq:direct_interp}
	P_{ij} = - A_{ij} ~\frac{\sum_{k \in N_{i}} A_{ik}}{ A_{ii} \sum_{k \in C_{i}^{s}}A_{ik}}
\end{equation}
where $N_{i}$ denotes the neighbors of vertex $i$, and $C_{i}^{s}$ denotes the
strong, coarse neighbors of vertex $i$.  Specifically $j \in C_{i}^{s}$
if and only if $j\in N_i$, $j\in \mathcal{C}$ and $S_{ij} - \tau > 0$.
If $i$ is instead a C-vertex,
\begin{equation*}
	P_{ij} = \begin{cases}
             1 & j = i \\
             0 & j \neq i
           \end{cases} .
\end{equation*}
A derivation of \eqref{eq:direct_interp} can be found in~\cite{de_sterck_distance-two_2008}. One can see that the $i^{th}$ row of 
$P$ is simply a weighted sum where weights are proportional to $A_{ij}$ and normalized by the fraction term,
This fraction does not depend on $j$ and guarantees that the sum of the nonzeros in the $i^{th}$ row of $P$ 
$(~= \sum_{k \in C_{i}^{s}} P_{ik}~)$ is equal to one when the sum of $A$'s nonzeros in the $i^{th}$ row is zero.
This is generally true when $A$ represents differentiation as the derivative of a contant function is zero,
which is effectively approximated by multiplying $A$ by a vector where all entries are constant.

The GNN representation of this algorithm, which is shown in
Table~\ref{algo:direct_interp}, requires two layers. We use
the coarse/fine splitting information and the matrix diagonal as fixed vertex
features. We also need a mutable, temporary vertex feature $\alpha$ for communicating
information between the two GN layers. As edge features, we use the off-diagonal
entries of $A$ and the strength of connection information as fixed features, and
the output weight, $P_{ij}$ is edge-specific and mutable.

As the algorithm proceeds, the edge update in the first layer is used to pass
the coarse/fine splitting information,
\begin{equation*}
 C_j = \begin{cases}
   1 & j\in\mathcal{C}\\
   0 & j\in\mathcal{F}
   \end{cases},
\end{equation*}
of the receiving vertex to the edge. Next, the edge aggregation is used to perform
the necessary summations since all the quantities are specific to the $i$th
vertex. In the vertex update of the first layer, the $A_{ii}$ factor is multiplied
through. Finally, in the second layer, only the edge update function is
utilized, which incorporates the $A_{ij}$ value into the weight and also
assigns 
zeroes to all rows associated with coarse points. After the computation of the GNN is
complete, some additional post-processing is necessary. As mentioned previously,
the GNN generates an updated graph with the same structure as the input graph
(which is square since it is derived from $A$ and $S$), but the prolongation
operator is rectangular. To get a complete prolongation operator, two additional
operations would be necessary: (a) setting the diagonal to one and (b) removing
all columns associated with fine vertices. The first of these operations can be
wrapped into the edge update in the second layer of the GNN network, but the
second must be completed after the computation of the entire graph neural
network. In doing so, the standard direct interpolation operator is obtained.

\begin{table}
  \caption{Direct Interpolation Kernel as a Graph Network}\label{algo:direct_interp}
  \begin{gnntable}
    \GNNDataAndLayer
    {$A_{ij}, \hat{S}_{ij}$}{$v_{ij}, P_{ij}$ [Output]}
    {$A_{ii}, C_i$}{$\alpha_i$}
    {---}{---}
    {Layer 1}
    {$v_{ij} = C_{j}$} {$\alpha_i = \frac{1}{A_{ii}} \overline{\gamma}_{i}$}{---}
    {$ \overline{\gamma}_{i} = \phantomLeft \displaystyle \frac{\sum_j A_{ij}}{\sum_j A_{ij} v_{ij} \hat{S}_{ij} } \phantomRight$}{---}{---}
    \GNNExtraLayer{Layer 2}
    {$P_{ij} = (1-C_i)(-A_{ij}\alpha_i)$}{---}{---}
    {---}{---}{---}
  \end{gnntable}
\end{table}

\section{Training Within GNN Frameworks}\label{sec:iter-methods}
The GNN models considered in the previous section used prescribed aggregation and update functions. 
Typically, update functions used are parameterized multi-layer perceptron networks, while the variadic aggregation functions
are prescribed. This section gives two examples that illustrate how learned parameters 
can be introduced into update functions so that the GNN model can be trained to perform a numerical task.

\paragraph{\textbf{Code Availability}}
The code for the examples is this chapter is provided in the repository at \url{https://github.com/sandialabs/gnn-applied-linear-algebra/}. Only PyTorch code is provided. The repository includes code to create all datasets, implement the given GNNs, and train the models.

\subsection{Learning a Diagonal for the Jacobi Iteration}\label{sec:jacobi_nn}
As discussed, an iteration of the weighted Jacobi method for solving $A x = b$ is defined by
\begin{equation} \label{eqn:Jacobi repeat}
x^{k+1} = x^{k} + \omega D^{-1}(b-A x^{k}),
\end{equation}
where the solution at the $k^{th}$ Jacobi iteration is updated based on the associated residual.
The $i^{th}$ entry of the residual is scaled by $ \omega / D_{ii} $ to update the
$i^{th}$ solution entry.  Here,  $\omega$ is a user-provided scaling value and $D$ is a diagonal matrix
whose nonzeros are given by $D_{ii} = A_{ii}$.  Ideally, $\omega$ is chosen so that
convergence is attained in the fewest number of iterations. While this ideal $\omega$ can be 
determined by computing eigenvalues of the matrix $D^{-1} A$ when $A$ is a symmetric positive definite matrix, 
in practice it is often chosen in an ad hoc fashion.

We pursue an alternative machine learning approach that selects a diagonal relaxation operator $\bar{D}$ based on the matrix.
A generalized Jacobi method (or scaled Richard iteration) is given by 
\begin{equation} \label{eqn:generalized Jacobi}
x^{k+1} = x^{k} + \bar{D}^{-1}(b-A x^{k}),
\end{equation} 
where the diagonal matrix $\bar{D}$ is not defined by $A$'s diagonal entries but is instead
chosen by a machine learning algorithm to reduce the number of iterations required to reach convergence.
We follow the framework introduced in~\cite{taghibakhshietal2022_LearningInterfaceConditions} for learning the diagonal
of the generalized Jacobi iterative method, though we consider a different target class of matrices to demonstrate
the approach focusing on Jacobi as a relaxation method in multigrid. Within multigrid, the relaxation method should reduce high frequency errors
as a complement to the coarse grid correction that addresses low frequency errors. 

\begin{figure}[ht]
  \centering
  \caption{An example mesh with a thin band of elements used for training and testing. Here $h = \frac{1}{7}$ and $\beta = 0.05$}
  \includegraphics[width=0.7\textwidth]{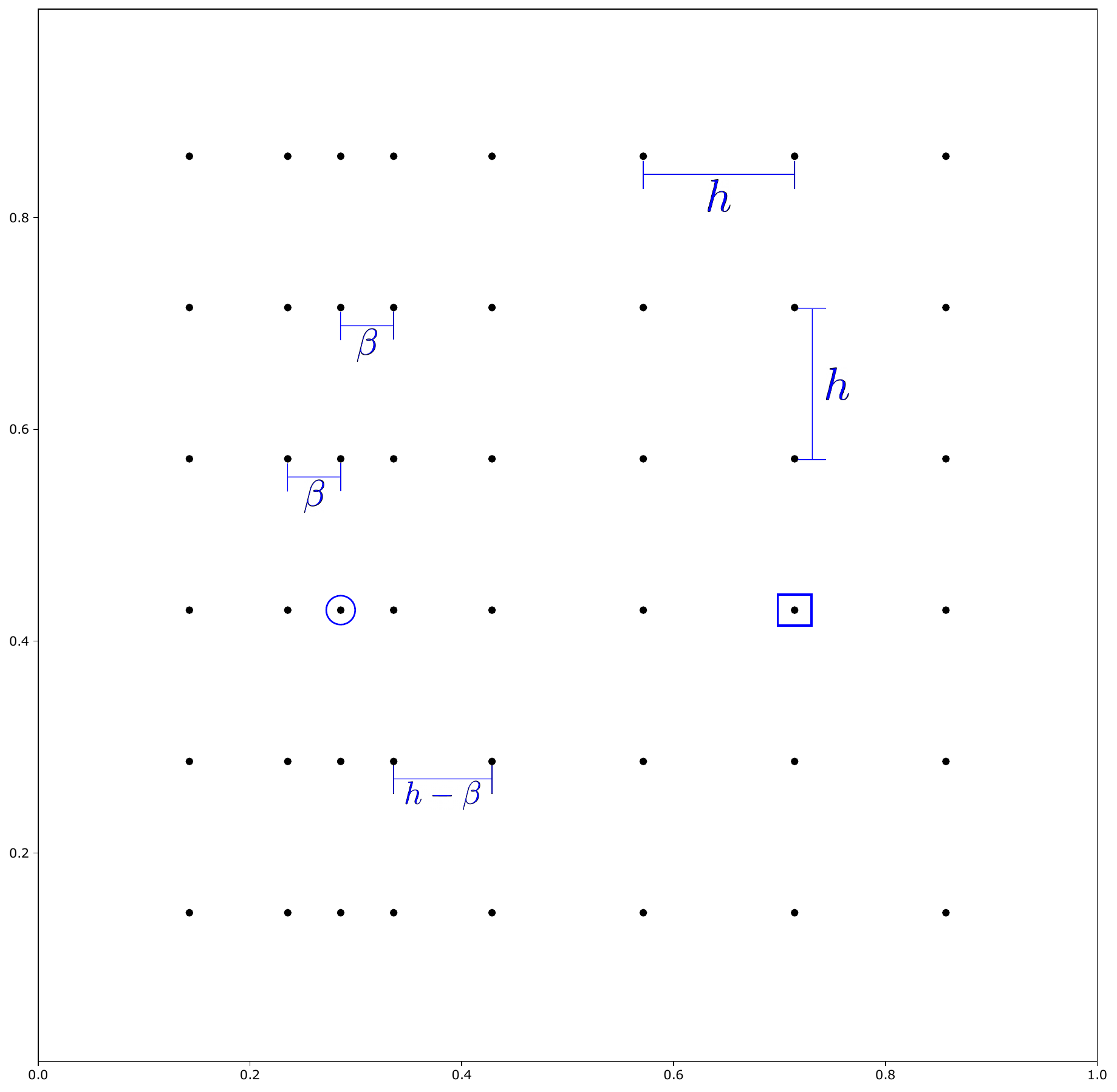}
  \label{fig:JacobiMesh}
\end{figure}

As an example, we discretize a 2D Poisson operator
\begin{equation*}
  \frac{\partial^2 u}{\partial x^2} + \frac{\partial^2 u}{\partial y^2 }
\end{equation*}
on the unit square $[0,1]\times[0,1]$ with homogeneous Dirichlet boundary conditions: $u(x,0) =
u(0, y) = u(x,1) = u(1, y) = 0$. The domain is tiled using linear finite elements with uniformly shaped quadrilaterals,
except for a small, thin band of tall and skinny quadrilaterals,
such as the mesh shown in Figure~\ref{fig:JacobiMesh}. Note that the horizontal location of the
band and its width can be changed to produce different variations of the problem. This
example is meant to highlight a situation where a generalized $\bar{D}$ matrix
within Jacobi's method might be advantageous, and should not be misconstrued as good meshing practice.  

In the matrix arising from discretizing this system, the abrupt change in mesh spacing yields
significant differences when comparing rows associated with mesh nodes far from the band
versus mesh nodes within the band. Two matrix stencils that highlight this variation are 
$$
\frac{1}{6 } \StencilTwoD{16}{-2}{-2}{-2} \hskip .5in \mbox{and} \hskip .5in 
\frac{\beta}{ 6h }\StencilTwoD{\phantom{+}8 + 8 \left( \frac{h}{\beta}\right)^2}{-4 + 2 \left(\frac{h}{\beta}\right)^2} {\phantom{+}2 - 4\left(\frac{h}{\beta}\right)^2}{-1-\left(\frac{h}{\beta}\right)^2} .
$$
The left stencil corresponds to a mesh node like the node surrounded by a square in 
Figure~\ref{fig:JacobiMesh} which is away from the boundary and the thin band (so locally the mesh is uniformly spaced in both
the $x$ and $y$ directions).
The right stencil corresponds to a
mesh node like the node surrounded by a circle in Figure~\ref{fig:JacobiMesh} which in the thin 
band and away from the boundary; the vertical spacing is $h$ and the
horizontal spacing (both to the left and right of the node) is $\beta$.
When $h/\beta = 1$, the left and right stencils correspond.
Notice that for $\beta \ll h$, the $h/\beta$ terms dominate and one can see that the 
left and right stencils  
are very different. 
While our example is artificially devised to stress Jacobi's
method, this type of abrupt stencil change occurs in more realistic scenarios. For instance in~\cite{DOMINO2018331},
generalized Jacobi iterations are used for a discontinuous Galerkin discretization where penalties
glue meshes together. 

\subsubsection{Datasets}\label{sec:jacobi:dataset}
  The training, test and validation datasets are
  generated by varying the location and width of the band in
  the mesh. The steps for producing the datasets are as follows: 
  \begin{enumerate}
    \item Create a 2D uniform mesh which contains $N_{y}$ points in each direction
    \item Randomly select $\beta$ between $\beta_{\text{min}}$  and $\beta_{\text{max}}$ 
    \item Randomly select an existing $x$-coord from the mesh in step 1
    \item Build the band mesh by placing new mesh vertices a distance of $\beta$ to the left and right
    of all vertices with the $x$-coordinate in step 3.
    \item Build the matrix based on a finite element discretization
    \item Repeat steps 1-5 to generate 1000 matrices
    \item Matrices 1-800 become the training dataset
    \item Matrices 801-850 become the validation set
    \item Matrices 851-1000 become the test set
  \end{enumerate}
  The training set defines the loss function used by gradient descent to determine model parameters 
  that yield a reduced loss value. 
  Using the validation set, the model (built
  with the training set) can be evaluated over a range of different model and algorithm choices
  (referred to as hyperparameters) such as model architectures and optimizers. The test dataset is
  used to evaluate the final model obtained after all hyperparameter tuning has occurred.  The test
  set is never used during the model selection process.

  \subsubsection{Input/Output Attributes}
  The input edge attributes are the nonzero matrix entries $A_{ij}$ excluding the diagonal. 
  The matrix diagonal, $A_{ii}$ is included as a a vertex attribute for each $i$.
  The output vertex attribute will be $d_i = \bar{D}^{-1}_{ii}$ for the generalized Jacobi method. 
  See the \emph{Data} section of Table~\ref{algo:jacobi_nn}.
  
  \subsubsection{Loss Function}
  A neural network is trained by identifying parameter values that approximately minimize or at
  least significantly reduce the value of a loss (also known as an objective) function. 
  One standard machine learning technique involves finding {\it
  optimal} answers (typically from an expensive process) and then training the neural network to
  match these optimal values; a technique referred to as supervised learning. In our case, this
  would require computing a set of optimal $\bar{D}_{ii}$ values for each training scenario using
  a (possibly expensive) numerical procedure.  Then, a loss function would be defined that
  minimizes the difference between these pre-computed values and the GNN output version of these
  values. Unfortunately, a tractable numerical procedure for pre-computing the $\bar{D}_{ii}$ values
  is not apparent given the fact that numerous eigenvalue calculations would be needed to
  determine an {\it optimal} high dimensional vector (in $\mathbb{R}^{N}$) to define $\bar{D}$ for
  each training case.  
  
  Instead we define the
  loss function based on the desired performance of the Jacobi method itself. We would like
  to minimize the number of required Jacobi iterations to reach some specified convergence
  tolerance. However, efficient optimization methods for neural network training require
  computing the gradient of the loss function and since the number of iterations is a
  discrete quantity, we cannot calculate a gradient with respect to it. 
  Thus, the number of iterations cannot be used and a differentiable quantity
  that measures the performance of a Jacobi iteration will be used instead. From the theory of iterative methods,
  we that the asymptotic convergence rate of a linear iterative method is given by the spectral radius
  of the error propagation matrix. The error propagation matrix for generalized Jacobi is 
$I \hskip -.02in - \hskip -.02in \bar{D}^{-1} A $, which is obtained by substituting $e^k = x^k
\hskip -.02in - \hskip -.02in A^{-1} b$ and $e^{k+1} = x^{k+1} \hskip -.02in - \hskip -.02in A^{-1}
b$ in \eqref{eqn:generalized Jacobi}.
  Thus, we seek to
  minimize the spectral radius of this $I \hskip -.02in - \hskip -.02in \bar{D}^{-1} A $ matrix. 
  
  Here, we face another obstacle: standard automatic differentiation tools do not support eigenvalue
  algorithms. Even if support was added,
  \cite{wangetal2019_BackpropagationFriendlyEigendecomposition} shows that standard eigenvalue
  decomposition algorithms and power method algorithms give numerically unstable gradients. Thus, we
  need a different method to approximate the spectral radius of a matrix which gives a numerically
  stable eigenvalue estimate. In \cite{taghibakhshietal2022_LearningInterfaceConditions}, it is
  proved that the spectral radius of a matrix, $B$ can be approximated by 
  \begin{equation*}
    \max_{i=1, 2, \ldots , m} \{ \lVert B^{K} \hat{u}_i \rVert^{\frac{1}{K}} \} 
  \end{equation*}
  where $K$ is a user defined number of iterations, $\{\hat{u}_{i}\}_{i=1}^{m}$ is a set of vectors randomly chosen from 
  the surface of the unit sphere. 
  For multigrid relaxation, we instead want to maximize the 
  performance of Jacobi with regard to only high-frequency errors. 
  We can characterize the high frequency space using a discrete sine transformation matrix $V$ that is
  defined such that each column has the general form
  $$
      \alpha \sin (\theta_{x} \pi x) \sin (\theta_{y} \pi y) .
  $$
  Specifically, each of the $N_{y}^2$ columns corresponds to a unique $(\theta_{x},\theta_{y})$ pair
  chosen from $\theta_{x} = 1, ..., N_{y} $ and  $\theta_{y} = 1, ..., N_{y} $.  The scalar $\alpha$
  is picked to ensure that the norm of each column is one and for the results in Section~\ref{sec:TrainJacobi:Results}, we use $N_{y} = 38$.  Partitioning $V$ into low frequency and
  high frequency columns, we have $V = \left [ V_{lf} ~~ V_{hf} \right ]$ where the low frequency
  modes correspond to pairs where both $\theta_{x} \le N_{y}/2$ and $\theta_{y} \le N_{y}/2$ while
  all other modes define $V_{hf}$. As we seek to minimize high frequency errors via Jacobi
  relaxation, for each matrix $A^{(j)}$ we train the GNN to find $\bar{D}$ that minimizes 
  \begin{equation}
    \max_{i= 1, 2, \ldots, m} \{ \lVert (I - \bar{D}^{-1} A^{(j)} )^{K}  \hat{u}_{i}  \rVert^{\frac{1}{K}} \} %
  \end{equation}
  where each $\hat{u}_i$ is a randomly chosen column of $V_{hf}$.
  Therefore, if $N$ is the number of matrices in the dataset being evaluated (training, validation, or test), the overall loss function is
  \begin{equation}
    \mathcal{L} = \sum_{j=1}^{N} \max_{i= 1, 2, \ldots, m} \{ \lVert (I - \bar{D}^{-1} A^{(j)} )^{K}  \hat{u}_{i}  \rVert^{\frac{1}{K}} \} %
    \label{eq:TrainableJacobiLoss}
  \end{equation}
  For the results given in Section~\ref{sec:TrainJacobi:Results}, we use $K = 3$ and $m = 20$.

\subsubsection{Architecture}
  In general
  the choice of architecture requires a search over range of parameters known as an ablation study. 
  However, for this example, a fixed architecture for the neural network is chosen as a simple demonstration. 
  The GNN consists of a
  single GNN layer with a vertex update, as can be seen in
  Table~\ref{algo:jacobi_nn}.
  \begin{table}
    \caption{Training a Diagonal for Jacobi as a Graph Network}\label{algo:jacobi_nn}
    \begin{gnntable}
      \GNNDataAndLayer
      {$A_{ij}$}{---}
      {$A_{ii}$}{$d_i$ [Output]}
      {---}{---}
      {Layer 1}
      {---}{$d_i = \neuralnet(A_{ii}, \overline{c})$}{---}
      {$\overline{c} = [min, mean, sum, max]$}{---}{---}
    \end{gnntable}
  \end{table}
  The edge to vertex  aggregation function, $\rho_{e\rightarrow v}$, calculates the min, mean, sum and max of the edge
  attributes, resulting in four quantities per vertex.
  These four attributes are concatenated
  with the single vertex attribute ($A_{ii}$) to form a vector in $\mathbb{R}^5$ input into the vertex
  update function. The vertex update function, $\phi_{v}$ is a feed-forward neural network
  consisting of 3 layers. For more information about the architecture of this GNN, see Appendix~\ref{app:arch:jacobi}.

  \subsubsection{Training Methodology}
  The model is implemented in PyTorch~\cite{PyTorch2015}, using the PyTorch Geometric
  library~\cite{PytorchGeometric2019}. The GNN is trained for 62 epochs using the Adam optimization
  method~\cite{Adam2014} with batch size of 100. The number of epochs can have a large influence on
  the overall success of the training phase. Too few epochs and the training may not reach peak performance,
  while too many epochs could lead to the model "over-fitting". In this case, the model fits the data
  too tightly and so the ability of the model to generalize to new data degrades. 
  There are different methods to address over-fitting; we choose to stop training when the performance of the
  model on the validation set starts to degrade. 
  For this experiment,
  the number of epochs is selected by training the model for 100 epochs, but stopping the training
  when the loss on the validation set is minimized. By stopping when the validation error is
  minimized, we obtain good generalization properties for future data.

  \subsubsection{Results}\label{sec:TrainJacobi:Results}
\begin{figure}[htbp]
    \centering
    \caption{Top 10 largest eigenvalues of the generalized Jacobi
      relaxation method for two different test matrices,
      comparing the different methods of selecting the $\omega$
      parameter in weighted Jacobi to the learned diagonal from the
      GNN.}
    \subfigure{\includegraphics[width=0.4\textwidth]{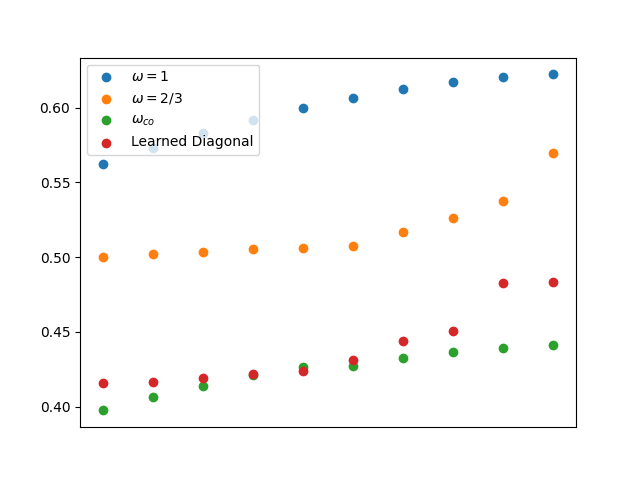}}
    \subfigure{\includegraphics[width=0.4\textwidth]{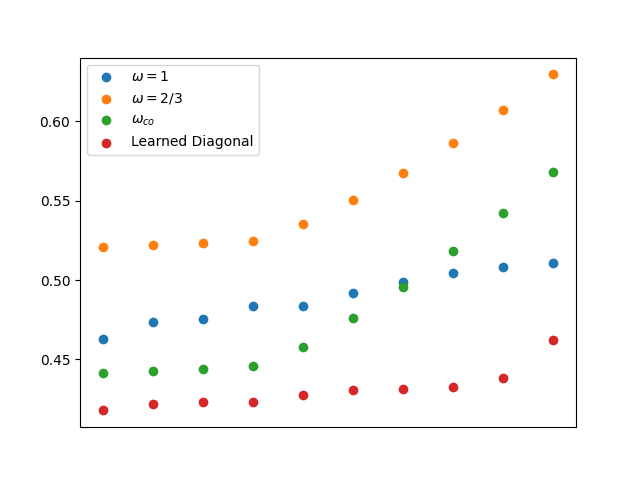}}
    \label{fig:Top10Eigs}
  \end{figure}

\begin{figure}[htbp]
    \centering
    \caption{Distributions of the differences in the maximum eigenvalues from the learned diagonal and the given constant weights for the Jacobi method. Counts left of the 0.0 line indicate matrices for which the associated method has lower maximum eigenvalue than the learned method (and by how much), while counts right of the 0.0 line indicate matrices for which the associated method has larger maximum eigenvalue than the learned method.}
    \includegraphics[width=0.9\textwidth]{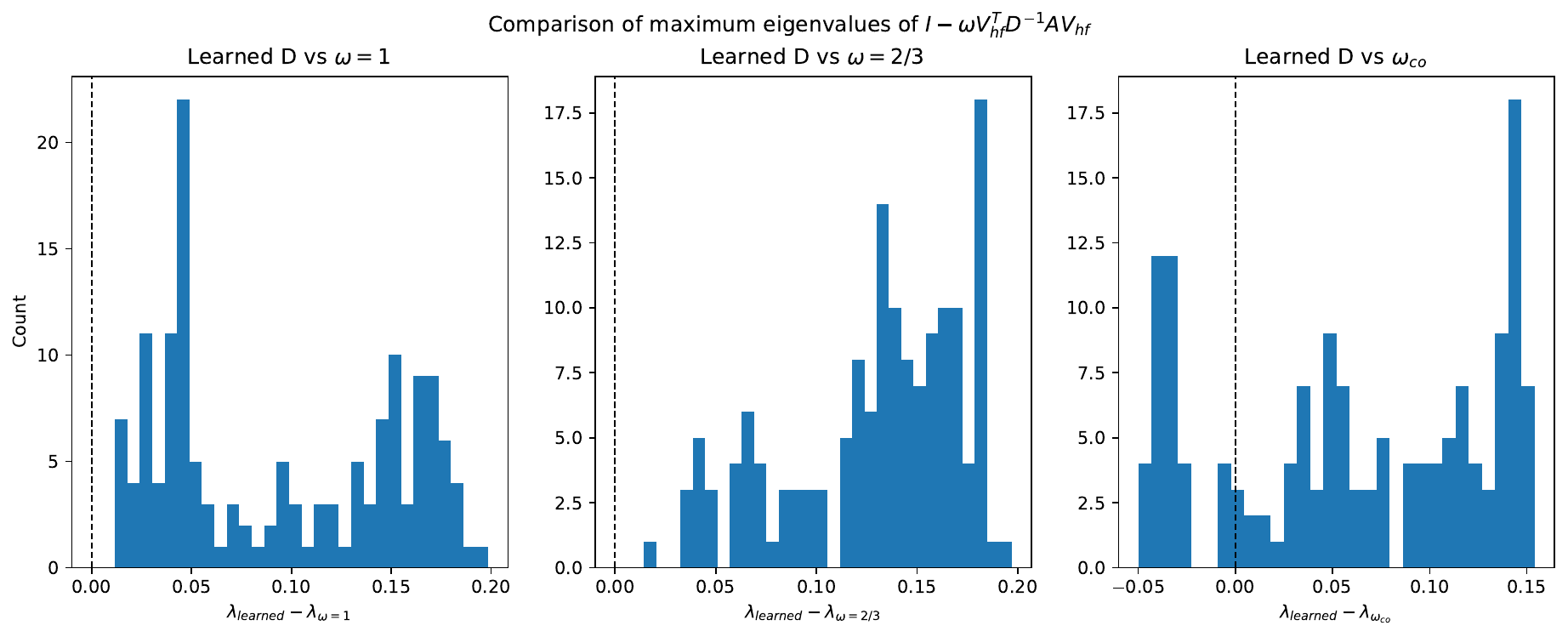}
    \label{fig:JacobiComparison}
  \end{figure}

  We compare the result of the generalized Jacobi iteration using the learned model $\bar{D}$ with a
  traditional weighted Jacobi using three different standard weights: standard Jacobi ($\omega =
  1$), a ``default'' weight ($\omega = \frac{2}{3}$), and a ``classical optimal'' weight
  ($\omega_{co}  = \frac{2}{\lambda_{\text{min}} + \lambda_{\text{max}}}$ where
  $\lambda_{\text{min}}, \lambda_{\text{max}}$ are the minimum and maximum value of $D^{-1} A$
  respectively). Note that the ``classical optimal'' weight is only optimal when Jacobi is used
  directly to solve $A x = b$ as opposed to Jacobi relaxation within a multigrid algorithm. To
  compare the different Jacobi relaxation procedures, we compute the maximum eigenvalue(s) of
  $I\hskip -.03in - \hskip -.03in\omega V_{hf}^{T}D^{-1} A V_{hf}$  (or of $I \hskip -.03in - \hskip
  -.03in V_{hf}^{T}\bar{D}^{-1} A V_{hf}$ for the generalized Jacobi method), thereby investigating
  each method's ability to damp high frequency error. The 10 largest eigenvalues associated with two
  of the test matrices are shown in Figure \ref{fig:Top10Eigs}.  
  We can see from the two sample matrices that in some cases the learned diagonal outperforms the
  rest of the methods, while in some cases, the $\omega_{co} $ weight performs the best. It would
  therefore be helpful to see a full comparison across all test matrices. In
  Figure~\ref{fig:JacobiComparison},
  we compare the distributions of maximal eigenvalues of the standard methods with respect to the
  learned diagonal method. Hence, in the figure, all counts left of the 0.0 line indicate matrices
  where the associated method performed better than the learned method, while all counts to the right
  of the 0.0 line indicate matrices where the associated method performed worse than the learned
  method. From these results, we see that the learned method gives a better result than the standard
  non-weighted and $\omega = \frac{2}{3}$ weighted Jacobi method in all cases and outperforms the
  $\omega_{co}$ weight in more than 75\% of the test matrices. 
  
\begin{figure}[htbp]
    \centering
    \caption{A plot of all the test matrices according to the width and location of their band, where color represents the winning method. Red represents the learned diagonal method while green represents the $\omega_{co}$ constant weight method.}
    \includegraphics[width=0.5\textwidth]{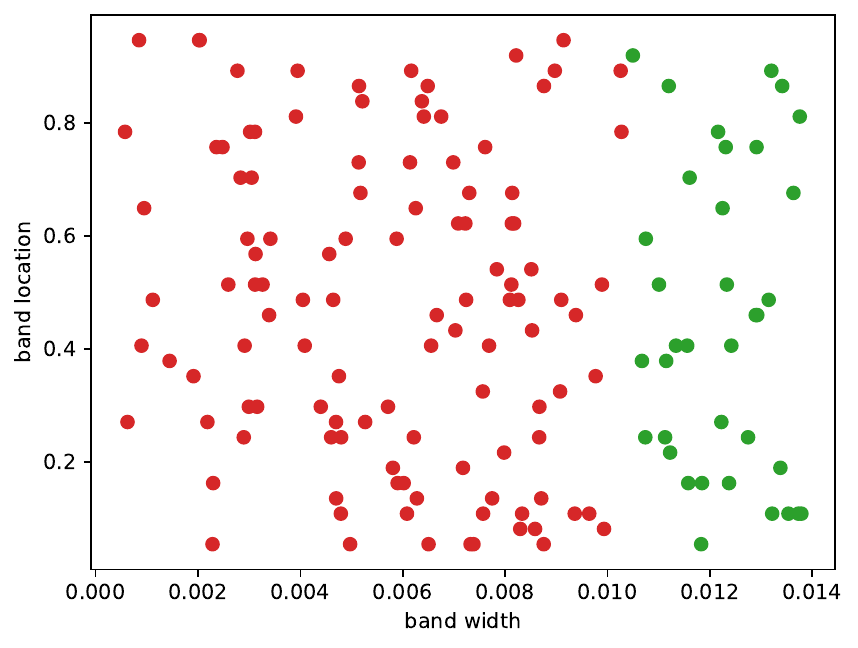}
    \label{fig:JacobiWinners}
  \end{figure}

  The shape of the right-most histogram
  in Figure~\ref{fig:JacobiComparison} suggests that there might be areas of the test space where
  the learned method performs better and areas where the $\omega_{co}$ performs better. To discover
  what patterns might exist, we map test matrices with regard to the width and location of the band 
  in Figure~\ref{fig:JacobiWinners}
  and color their
  associated dot with the top performing method. We can see from this plot that there is a clear
  pattern with regard to which matrices are more effectively solved using the learned diagonal
  versus the $\omega_{co}$, constant weight. This makes intuitive sense because for larger widths,
  the mesh is ``more uniform'' and thus the correctly chosen constant will perform well, while for
  smaller widths, allowing the weight to vary on a per-row basis allows for more flexibility to address
  the more drastically different elements.

\begin{figure}[htbp]
    \centering
    \caption{Diagonals at each mesh point for the $\omega_{co}$ method and the learned method}
    \subfigure{\includegraphics[width=0.48\textwidth]{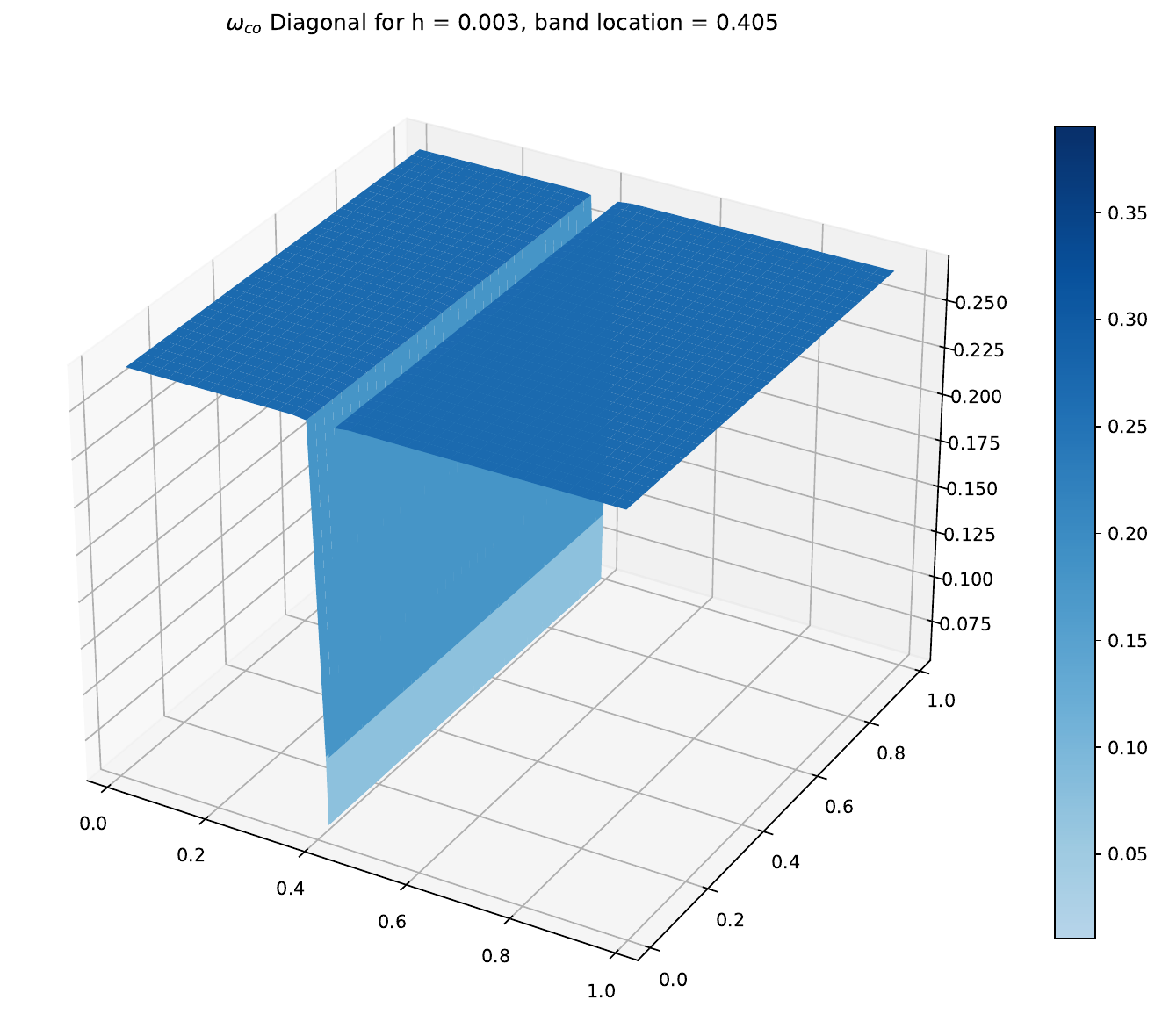}}
    \subfigure{\includegraphics[width=0.48\textwidth]{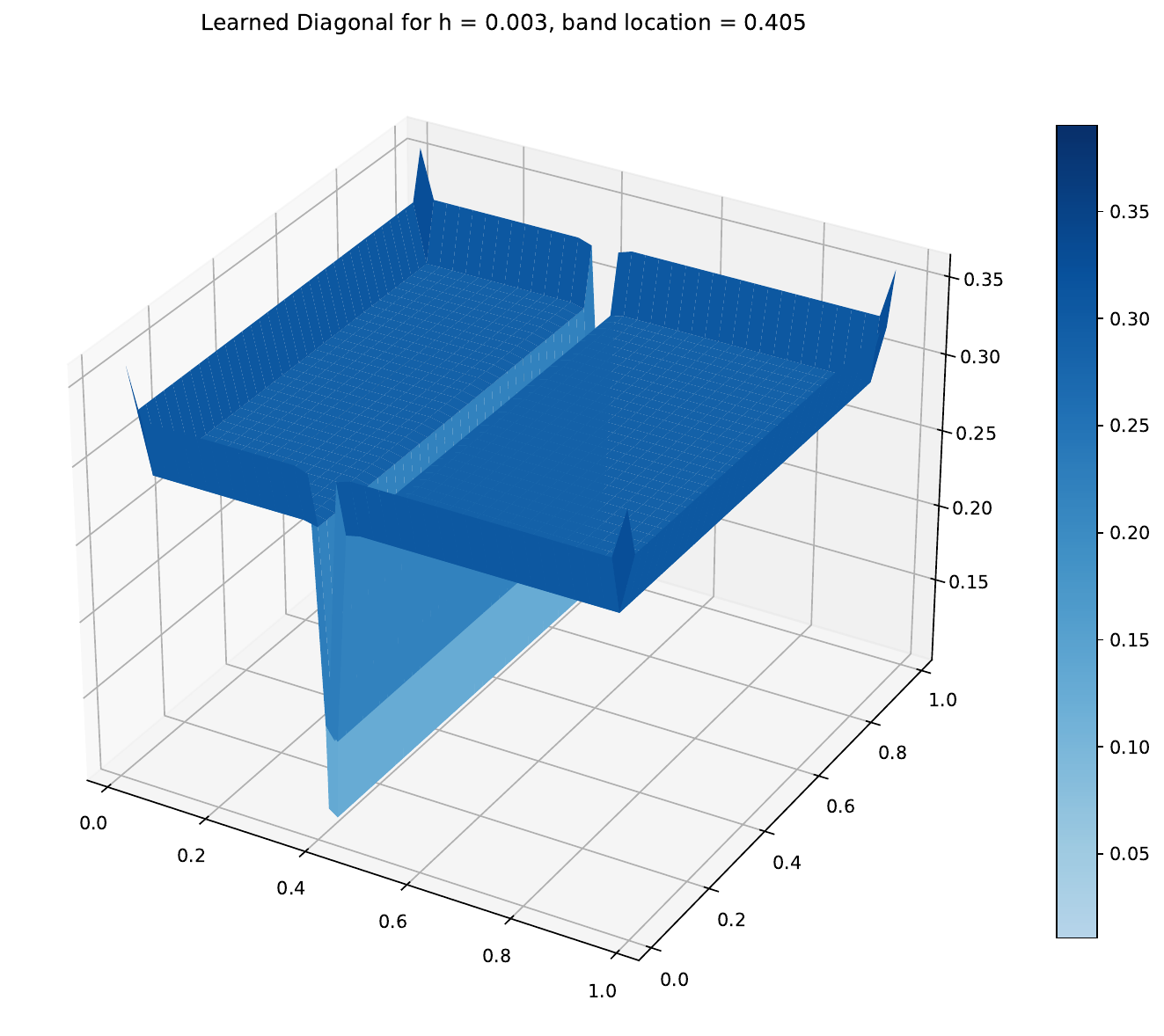}}
    \label{fig:JacobiDiagSurface}
  \end{figure}

  Finally, comparison plots for the scaled diagonal values arising from the $\omega_{co}$ method and
  the learned diagonal are given in Figure~\ref{fig:JacobiDiagSurface} for a test matrix with $h =
  0.003$ and where the band is approximately located at $x=0.405$.
  There are several observations which can be made from the images. First, we observe that in
  general, the learned method yields larger diagonal values than the $\omega_{co}$ method. It can
  also be observed that the learned method produces a larger diagonal in the center of the band than
  the $\omega_{co}$ method. Most drastically, we observe in the learned method that boundary points
  (especially the corners) have larger diagonal values than the interior points. This is true
  even in the band, where the boundary points are also larger than the interior points. This is a
  modification that the traditional weighted Jacobi method cannot replicate since the diagonal of
  $A$ is the same at the boundaries and interior.

  The shapes in Figure~\ref{fig:JacobiDiagSurface} demonstrate why the learned method is able to
  outperform the ``classical optimal'' weighting for Jacobi's method. For the standard weighted
  Jacobi, only a constant multiple of $D^{-1} A$ is allowed, which forces the shape of the diagonal
  to match the shape of the inverse of the diagonal of $A$. This shape cannot be adjusted, nor can
  the relative distances be changed (only scaled). In our new learned method, no such limitation
  exists, hence the space from which the learned method can select the diagonal values is much
  broader. This broader space presents greater opportunity to find an optimal diagonal.

\subsection{Example: Training a GNN to Determine Diffusion Coefficients}\label{sec:diff_training}
The strength of connection algorithms discussed in Section~\ref{sec:algo_soc_sa} 
are heuristic schemes with several shortcomings. One could
instead consider GNN models that learn a more sophisticated measure of strength. In this section, we
consider the related task of determining diffusion coefficients given
matrix stencil and coordinate information. Though not equivalent to determining strength
measures, there is a rough connection between the relative sizes of diffusion coefficients in
different directions and the strength of the matrix connections aligned with these directions.

Consider the diffusion operator 
\begin{align}
  \label{eq:DiffPDE}
\sum_{i = 1}^2 \sum_{j=1}^2 \frac{\partial}{\partial x_i} 
  \left[ D_{ij}(x, y) \frac{\partial 
      }{\partial x_j} \right] 
\end{align}
defined on the unit square with 
periodic boundary conditions.
$D(u,x,y)$ is assumed to have the form
\begin{equation*}
	D(x, y) = \begin{bmatrix}
                 \alpha(x,y) & 0 \\
                 0 & \beta(x,y)
               \end{bmatrix} .
\end{equation*}
The goal of this learning task
is to predict the diffusion coefficients $\alpha(x,y)$ and $\beta(x,y)$ at each grid point using information from the coefficient
matrix, coordinates of vertices, and mesh spacing.  

\subsubsection{Datasets}\label{sec:jacobi:diffusion} 
All the datasets consist of data generated from discretizing \eqref{eq:DiffPDE} with the finite element method. 
Using the following six-step procedure, $1000$ matrices are generated:
\begin{enumerate}
  \item Select a random integer $N \in [80,100]$ and let the mesh resolution be $h = \frac{1}{N}$
  \item Select $\theta_{\alpha, x}, \theta_{\alpha, y}, \theta_{\beta, x}, \theta_{\beta, y}$ each with
  uniform probability from $\{i : i \in \mathbb{Z} \text{ and } 0 \leq i \leq 6\}$ 
  \item Define $\alpha (x,y) = \cos \left(\theta_{\alpha , x} \pi x \right)^2 \cos \left(
  \theta_{\alpha, y}  \pi y \right)^2 $ 
  \item Define $\beta (x,y) = \cos \left(\theta_{\beta , x} \pi x \right)^2 \cos \left(
  \theta_{\beta, y}  \pi y \right)^2 $ 
  \item Discretize \eqref{eq:DiffPDE} on a uniform 2D quadrilateral mesh with resolution $h$ in both $x$ and $y$
  directions to construct the matrix operator $A$
  \item Generate input features and output targets for this matrix
\end{enumerate}
Notice that diffusion coefficients are chosen so that there is no discontinuity over the periodic 
boundary conditions. Thus, the diffusion fields vary smoothly over the domain.

The data set generated is divided into
a training set (700 matrices), validation set (200 matrices), and
test set (100 matrices). The training set is used to determine model parameters that yield a
sufficiently small loss. The validation set is used to determine the model
architecture, as will be explained shortly. Finally, the test set is only used once at the end of the
study to evaluate the final model's effectiveness on unseen data. The test set is evaluated in
the results sections where it is used as a trustworthy indication of how the chosen model will
perform on new, completely unseen data.

\subsubsection{Input/Output Attributes}
The left column of Table~\ref{algo:diff_nn} describes the input/output attributes. 
The vertex attributes $v_{i}$ consist of the
matrix diagonal entries. That is, $v_{i} = A_{ii}$. The edge attributes, $c_{ij}$ are constructed as
$c_{ij} = (A_{ij} , x_{\text{rel}}, y_{\text{rel}})$ where $x_{\text{rel}}, y_{\text{rel}} $ are the
relative differences between vertex $i$ and vertex $j$ in the $x$ and $y$ directions respectively,
scaled by $\frac{1}{h}$. For example, if vertex $j$ is the vertex southeast of vertex $i$, then the
edge features for the edge from vertex $i$ to vertex $j$ are $c_{ij} = (A_{ij}, 1, -1)$. There is
a global attribute in this example, which is the mesh resolution: $g = h$. The output of
the neural network is the updated vertex attributes which are the predicted values for $\alpha$ and
$\beta$ at each mesh vertex.

\subsubsection{Loss Function}
The loss function is given by the mean squared error (MSE)
$$
    {\cal L}_\omega^{(k)} = \frac{1}{2 N^2} \sum_{i=1}^{N^2} || \alpha^{(k)}(x_i,y_i) -  \tilde{\alpha}_\omega^{(k)}(x_i,y_i) ||_2^2 + 
                                || \beta^{(k)}(x_i,y_i) -  \tilde{\beta}_\omega^{(k)}(x_i,y_i) ||_2^2
$$
where the $(k)$ superscript denotes the matrix from the training set; $\omega$ refers to the model's trainable parameters (determined by numerical optimization during training); 
the $(x_i, y_i)$ correspond to different mesh points
and $\tilde{\alpha}_\omega^{(k)}$ and $\tilde{\beta}_\omega^{(k)}$ represents the GNN model predictions using the model parameter $\omega$.

\subsubsection{Architecture} \label{sec:arch:diffusion} 
The model architecture employs an encoder, followed by the GNN outlined in Table~\ref{algo:diff_nn}.
The encoder is applied to the edge and vertex attributes separately before the execution of the graph
neural network. Specifically,
an MLP, $\mathrm{Encoder}_e$, is applied to each set
of edge attributes. Similarly, an MLP, $\mathrm{Encoder}_v$, is applied to each set of vertex
attributes. Finally, a third MLP, $\mathrm{Encoder}_g$ is applied to the global attributes.  Details on
the design of the encoder MLPs can be found in Appendix~\ref{app:arch:diffusion}. 
Generically, encoders and decoders are common in GNNs. These are used to enhance the expressiveness of
 user input attributes with the aim to improve predictions and network trainability . 
Within the graph neural network layer, MLP neural networks are used for the edge update  and
 vertex update functions.

\begin{table}
  \caption{Training a Graph Network to Determine Diffusion Coefficients}\label{algo:diff_nn}
  \begin{gnntable}
    \GNNDataAndLayer
    {---}{$c_{ij}$}
    {---}{$v_i$ [Output]}
    {---}{$g$}
    {Encoder}
    {$c_{ij} = \mathrm{Encoder}_e(c_{ij})$}{$v_i ~= \mathrm{Encoder}_v(v_i)$}{$g ~~= \mathrm{Encoder}_g(g)$}
    {---}{---}{---}
    \GNNExtraLayer{Graph Neural Network}
    {$c_{ij} = \subscriptneuralnet{e}(c_{ij}, v_{i}, v_{j}, g)$}{$v_i ~= \subscriptneuralnet{v}(v_{i}, \overline{c}_i,g)$}{---}
    {$\overline{c}_i = [min, mean, sum, max]$}{---}{---}
  \end{gnntable}
\end{table}

This architecture was chosen after performing 
experiments using the
validation set.  The validation set allows researchers to experiment with different model 
parameters, architectures, optimizers, etc. without exposing the model to the test set prematurely. 
In our case, we 
investigated the following 
variations
\begin{itemize}
  \item Number of GNN layers: 1, 2, or 3
  \item Number of MLP layers in the update functions: 1, 2, 3, or 4
  \item Width of the hidden layers in the MLP update functions: 8, 16, 32, or 64
  \item Encoder architecture: No encoder, 1 layer with 16 neurons, or 3 layers with 16 neurons
  \item Decoder architecture: No decoder, 1 layer with 16 neurons, or 3 layers with 16 neurons.
\end{itemize}
This yields a total of $3\cdot4\cdot4\cdot3\cdot3 = 432$ total architectures combinations. 
To improve computational
efficiency, we evaluated the different models in two stages. In the first stage, we employ 25\%
of the training set and 25\% of the validation set to quickly identify the highest performing
models among the 432 possibilities. The five best-performing models are then evaluated in stage
two on the full validation set after they have been optimized using the full training set.
As described in more detail in Appendix~\ref{app:arch:diffusion}, the
best performing model of the 432  
is parameterized by $14,002$
trainable parameters (1 GNN layer, 2 MLP layers of width
16 in the update function and 3 layer 16 neuron encoders and decoders).

\subsubsection{Training Methodology}
The model is again implemented in PyTorch \cite{PyTorch2015}, using the PyTorch Geometric library
\cite{PytorchGeometric2019}. 
The Adam optimization method
\cite{Adam2014} is utilized to train the neural network weights with a batch size of $10$. 
Of particular importance is the choice of the  number of training epochs. As noted earlier, too few epochs 
may lead to poor accuracy while 
too many epochs may lead to over-fitting where the model performs 
poorly on
new, unseen data. 
We again choose to stop training when the performance of the
model on the validation set starts to degrade. 
In our case, the minimum validation loss occurs after
training the model for 187 epochs.

\subsubsection{Results}

\begin{figure}[ht]
  \centering
  \caption{Loss plots for training diffusion coefficients}
  \subfigure[Training loss for training diffusion coefficients]{\includegraphics[width=0.45\textwidth]{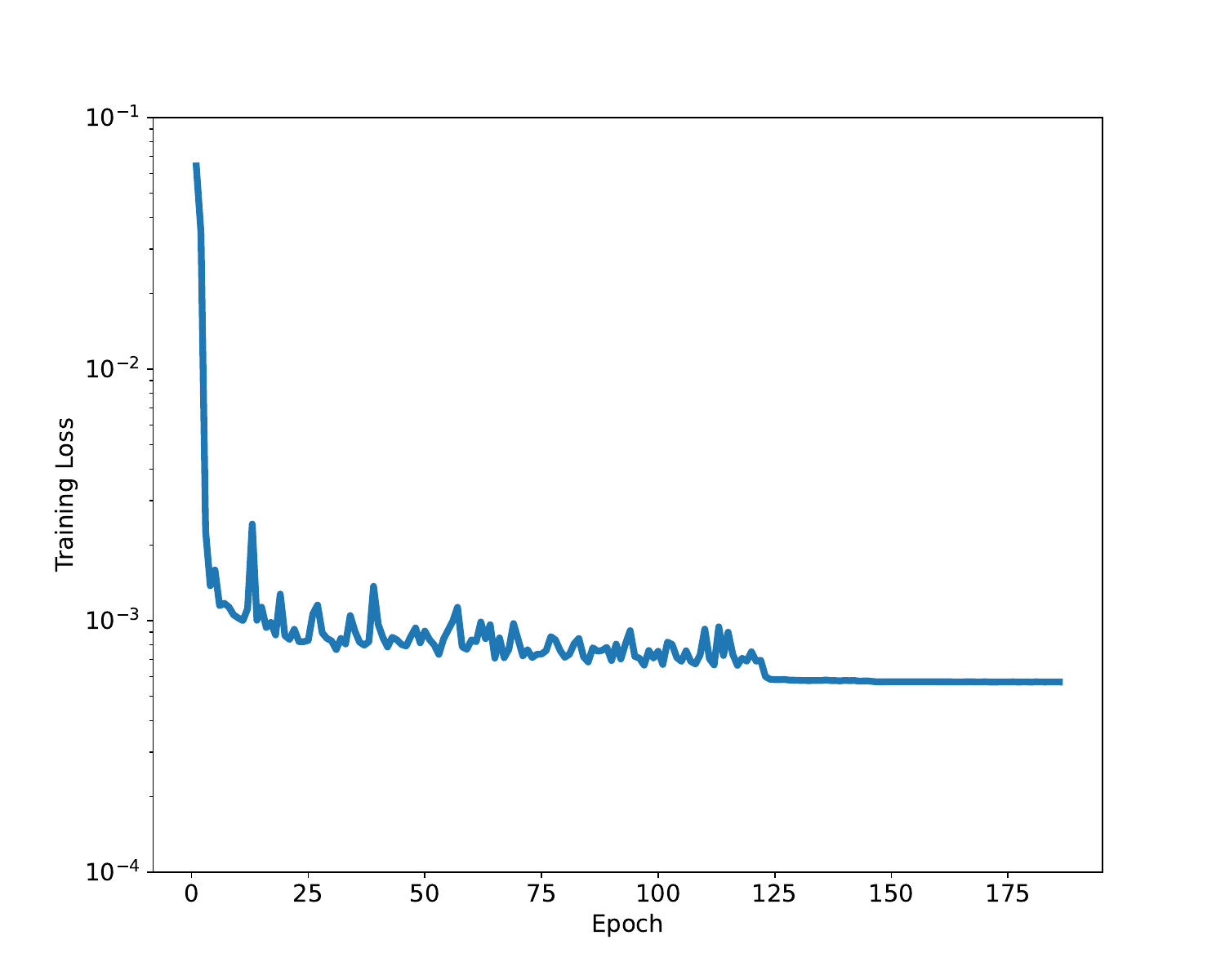}}
  \subfigure[Validation loss for training diffusion coefficients]{\includegraphics[width=0.45\textwidth]{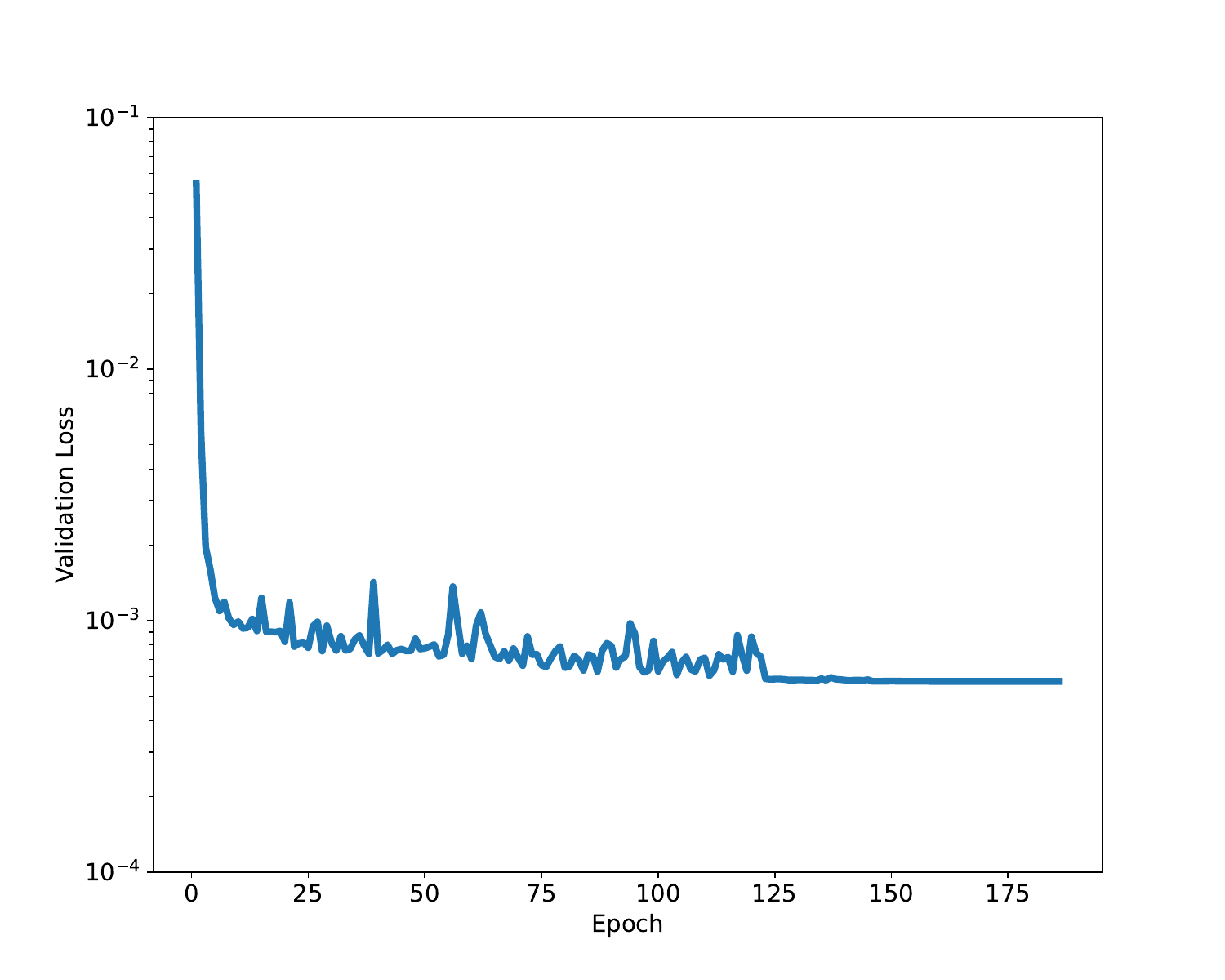}}
  \label{fig:training_diff_coeffs}
\end{figure}

Figure~\ref{fig:training_diff_coeffs} shows the MSE loss for the training and validation sets as a function
of number of epochs. Each epoch is a forward and backward propagation to compute the gradient over all batches in the
training set. Thus with a training set of $700$ entries, and a batch size of $10$ that is $70$ gradient computations
per epoch.
The loss
for both sets is attains a similar level ($5.73\times 10^{-4}$ for validation and $5.74\times 10^{-4}$ for testing).
This similarity indicates that it is unlikely that over-training has occurred. Over-training is usually marked by a descending training
error while simultaneously increasing validation error. This indicates the network does not generalize to data
outside of the training set.


\begin{figure}[ht]
  \centering
  \caption{Frequency vs mean-squared error for the trained model predicting point-wise diffusion
   coefficients from the finite element problem matrix and relative coordinates. The shaded region demonstrates the area where the training data was taken from, while the non-shaded region demonstrates where the model is extrapolating.}
  \includegraphics[height=4in,width=0.8\textwidth,trim={3.5cm 3cm 1.4cm 3.5cm},clip]{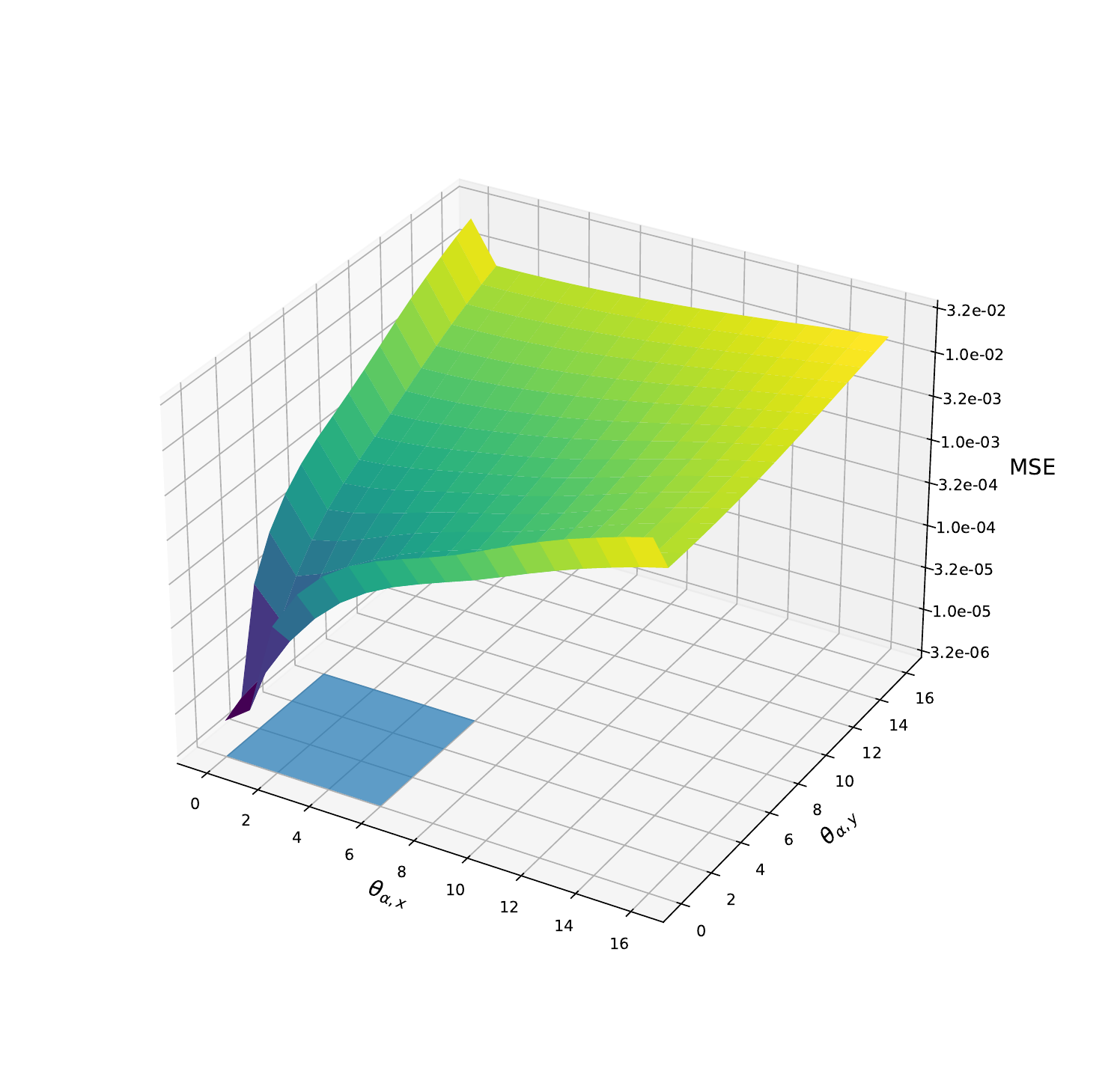}
  \label{fig:DiffFreqVsError}
\end{figure}

We also consider model performance in relation to the frequency
of the diffusion functions. To study this behavior, assume $D(x,y)$ takes the form
\begin{equation*}
  D(x,y) = \begin{bmatrix} \alpha(x,y) & 0 \\ 0 & \alpha(x,y) \end{bmatrix}
\end{equation*}
where $\alpha(x,y)$ has the form described in Section~\ref{sec:jacobi:diffusion}. Now, we can test the model
on problems where $\theta_{\alpha , x}$ and $\theta_{\alpha , y}$ are selected in a grid from
the set $\{i : i \in \mathbb{Z} \text{ and } 0 \leq i \leq 16\}$. The frequency versus
mean-squared error plot is given in Figure~\ref{fig:DiffFreqVsError}. Recall that the model is only
trained for $0 \leq i \leq 6$, so all frequency combinations outside this interval are being
extrapolated by the model. This portion of the subdomain is indicated by the shaded region.
All matrices used in this study are generated using the same finite
elements as previously, and all have 100 nodes in each direction yielding 10,000 by 10,000 matrix
problems. The figure indicates that coefficients with low frequencies, those in the training set, are well approximated. 
Deviating from the training set, there are two sources of potential error, both of which contribute to the increase in
the error observed in the figure. First, higher frequency coefficients suffer increasing error on a fixed ($100\times 100$)
mesh as the number of points per wavelength decreases. The second source of error comes from the model extrapolating
outside of the training set. In this context, examining the error away from the training frequencies we see that the error
remains relatively small even with larger departures from the training data.

Connecting this result to strength-of-connection metrics, we remark that 
in general, it can be difficult to correctly classify strong and weak connections  
 when $\alpha$ is small and $\beta$ is not small (or
vice versa) when linear finite elements are used on quadrilateral meshes.  Such a stencil where $\alpha(x,y) = 0.001$ and $\beta(x,y) = 0.8$ is shown below:
\begin{equation*}
  \StencilTwoD{1.068}{-0.533}{0.266}{-0.1335}
\end{equation*}  
Notice that if the following standard strength-of-connection metric is used:
\begin{equation*}
\abs{A_{ij}} \geq \theta \max_{k \neq i} \abs{A_{i,k}}  \Leftrightarrow A_{ij} \text{ is strong} 
\end{equation*}
then $\theta < 0.25$ would (incorrectly) classify all four cardinal directions as strong.
Allowing the trained model to predict the $\alpha$ and $\beta$ for this case, we obtain the
following values:
\begin{equation*}
  \alpha = -0.00052847 \quad \beta = 0.71802.
\end{equation*}
As such, the model has correctly predicted this stencil as having small $\alpha$ and with a
reasonable approximation for $\beta$. This provides an opportunity for giving better guidance
to direct the multigrid 
coarsening procedure.  Notice that the $\alpha$ prediction is slightly negative. This is a by-product of
using a LeakyReLU for the final activation allowing for negative outputs.
The standard ReLU, while strictly positive, yields several 
training cases that have smaller derivative values, slowing down learning, and ultimately
negatively impacting model accuracy.
Applying a thresholding similar to a ReLU as a post-training
step is possible, and would act as a more ``reasonable'' filter for outputs, but since we wished to
evaluate the performance of our model, no additional steps were applied to the output of the GNN.
This discrepancy also points to the possibility that more fine-tuning could improve the trained
model, but as this is a more educational paper in nature, we forego adding any additional
complexities to the model in this study.

\section{Conclusion}


We have examined graph neural networks (GNNs) from the
perspective of numerical linear algebra.  Our objective has been
to highlight the close relationship between several sparse matrix tasks
and GNNs. This relationship should not be surprising given the strong
connections between sparse matrices and graphs.  Several 
traditional numerical algorithms have been recast using a GNN representation.
The GNN algorithms presented here will often be much less efficient than  
specialized sparse linear algebra libraries like Trilinos~\cite{heroux2005overview}
or PETSc~\cite{balay2019petsc}, 
but are instead intended to familiarize the reader with the inner workings
and operations of GNNs by relating well-known algorithms to a
GNN counterpart. We believe that GNNs may prove useful for many 
sophisticated linear algebra tasks where traditional
algorithms have significant shortcomings. For example, one frequently recurring 
linear algebra theme is that of matrix approximation. The details
vary by context, but the general idea is to find an inexpensive surrogate matrix 
that approximates the behavior of a large sparse matrix.
For multigrid, this inexpensive matrix might be needed
to approximate the action of the original matrix on a low frequency space. In other 
situations (e.g., reduced order modeling) an approximation of matrix
vector products might be needed
only for vectors that lie in a relatively small subspace driven
by the simulation. These inexpensive approximations
leverage smaller, sparser or lower rank matrices.
Machine learning automates the task of generating approximation
using data in cases that can be difficult with
traditional approaches. 
The power of GNNs within machine
learning is that their flexibility can address general sparse matrices.
In this paper, we give 
two GNN learning examples to demonstrate, in simple cases, how GNNs
can address some linear algebra tasks that can be incorporated into more 
traditional methods.
While these examples are basic demonstrations of where GNNs can
be used, we hope they will provide inspiration for how to incorporate
GNNs into more complex problems.

\bibliographystyle{siamplain}
\bibliography{references,LearningInterp}

\begin{thebibliography}{10}

\bibitem{huggingfacetransformers}
{\em Hugging face transformers}.
\newblock \url{https://huggingface.co/docs/transformers/index}.
\newblock Accessed: 2023-08-14.

\bibitem{tensorflow2015-whitepaper}
{\sc M.~Abadi, A.~Agarwal, P.~Barham, E.~Brevdo, Z.~Chen, C.~Citro, G.~S.
  Corrado, A.~Davis, J.~Dean, M.~Devin, S.~Ghemawat, I.~Goodfellow, A.~Harp,
  G.~Irving, M.~Isard, Y.~Jia, R.~Jozefowicz, L.~Kaiser, M.~Kudlur,
  J.~Levenberg, D.~Man\'{e}, R.~Monga, S.~Moore, D.~Murray, C.~Olah,
  M.~Schuster, J.~Shlens, B.~Steiner, I.~Sutskever, K.~Talwar, P.~Tucker,
  V.~Vanhoucke, V.~Vasudevan, F.~Vi\'{e}gas, O.~Vinyals, P.~Warden,
  M.~Wattenberg, M.~Wicke, Y.~Yu, and X.~Zheng}, {\em {TensorFlow}: Large-scale
  machine learning on heterogeneous systems}, 2015,
  \url{https://www.tensorflow.org/}.
\newblock Software available from tensorflow.org.

\bibitem{balay2019petsc}
{\sc S.~Balay, S.~Abhyankar, M.~Adams, J.~Brown, P.~Brune, K.~Buschelman,
  L.~Dalcin, A.~Dener, V.~Eijkhout, W.~Gropp, et~al.}, {\em Petsc users
  manual},  (2019).

\bibitem{battaglia_relational_2018}
{\sc P.~W. Battaglia, J.~B. Hamrick, V.~Bapst, A.~Sanchez-Gonzalez,
  V.~Zambaldi, M.~Malinowski, A.~Tacchetti, D.~Raposo, A.~Santoro, R.~Faulkner,
  C.~Gulcehre, F.~Song, A.~Ballard, J.~Gilmer, G.~Dahl, A.~Vaswani, K.~Allen,
  C.~Nash, V.~Langston, C.~Dyer, N.~Heess, D.~Wierstra, P.~Kohli, M.~Botvinick,
  O.~Vinyals, Y.~Li, and R.~Pascanu}, {\em Relational inductive biases, deep
  learning, and graph networks}, arXiv:1806.01261 [cs, stat],  (2018),
  \url{http://arxiv.org/abs/1806.01261} (accessed 2021-05-19).
\newblock arXiv: 1806.01261.

\bibitem{bottou2018optimization}
{\sc L.~Bottou, F.~E. Curtis, and J.~Nocedal}, {\em Optimization methods for
  large-scale machine learning}, SIAM review, 60 (2018), pp.~223--311.

\bibitem{jax2018github}
{\sc J.~Bradbury, R.~Frostig, P.~Hawkins, M.~J. Johnson, C.~Leary,
  D.~Maclaurin, G.~Necula, A.~Paszke, J.~Vander{P}las, S.~Wanderman-{M}ilne,
  and Q.~Zhang}, {\em {JAX}: composable transformations of {P}ython+{N}um{P}y
  programs}, 2018, \url{http://github.com/google/jax}.

\bibitem{Briggs2000}
{\sc W.~Briggs, V.~E. Henson, and S.~McCormick}, {\em A Multigrid Tutorial},
  SIAM: Society for Industrial and Applied Mathematics, 2000.
\newblock Second Edition.

\bibitem{cyr2020robust}
{\sc E.~C. Cyr, M.~A. Gulian, R.~G. Patel, M.~Perego, and N.~A. Trask}, {\em
  Robust training and initialization of deep neural networks: An adaptive basis
  viewpoint}, in Mathematical and Scientific Machine Learning, PMLR, 2020,
  pp.~512--536.

\bibitem{de_sterck_distance-two_2008}
{\sc H.~De~Sterck, R.~D. Falgout, J.~W. Nolting, and U.~M. Yang}, {\em
  Distance-two interpolation for parallel algebraic multigrid}, Numerical
  Linear Algebra with Applications, 15 (2008), pp.~115--139,
  \url{https://doi.org/10.1002/nla.559},
  \url{http://doi.wiley.com/10.1002/nla.559} (accessed 2021-06-22).

\bibitem{DOMINO2018331}
{\sc S.~P. Domino}, {\em Design-order, non-conformal low-mach fluid algorithms
  using a hybrid cvfem/dg approach}, Journal of Computational Physics, 359
  (2018), pp.~331--351,
  \url{https://doi.org/https://doi.org/10.1016/j.jcp.2018.01.007},
  \url{https://www.sciencedirect.com/science/article/pii/S0021999118300172}.

\bibitem{AIElvis}
{\sc B.~Edwards}, {\em Hear elvis sing baby got back using ai --- and learn how
  it was made}.
\newblock
  \url{https://arstechnica.com/information-technology/2023/08/hear-elvis-sing-baby-got-back-using-ai-and-learn-how-it-was-made/}.
\newblock Accessed: 2023-08-14.

\bibitem{PytorchGeometric2019}
{\sc M.~Fey and J.~E. Lenssen}, {\em Fast graph representation learning with
  {PyTorch Geometric}}, in ICLR Workshop on Representation Learning on Graphs
  and Manifolds, 2019.

\bibitem{gilmer_neural_2017}
{\sc J.~Gilmer, S.~S. Schoenholz, P.~F. Riley, O.~Vinyals, and G.~E. Dahl},
  {\em Neural {Message} {Passing} for {Quantum} {Chemistry}}, arXiv:1704.01212
  [cs],  (2017), \url{http://arxiv.org/abs/1704.01212} (accessed 2021-08-16).
\newblock arXiv: 1704.01212.

\bibitem{glorot2010understanding}
{\sc X.~Glorot and Y.~Bengio}, {\em Understanding the difficulty of training
  deep feedforward neural networks}, in Proceedings of the thirteenth
  international conference on artificial intelligence and statistics, JMLR
  Workshop and Conference Proceedings, 2010, pp.~249--256.

\bibitem{goodfellow2016deep}
{\sc I.~Goodfellow, Y.~Bengio, and A.~Courville}, {\em Deep learning}, MIT
  press, 2016.

\bibitem{gori_new_2005}
{\sc M.~Gori, G.~Monfardini, and F.~Scarselli}, {\em A new model for learning
  in graph domains}, in Proceedings. 2005 {IEEE} {International} {Joint}
  {Conference} on {Neural} {Networks}, 2005., vol.~2, July 2005, pp.~729--734
  vol. 2, \url{https://doi.org/10.1109/IJCNN.2005.1555942}.
\newblock ISSN: 2161-4407.

\bibitem{he2015delving}
{\sc K.~He, X.~Zhang, S.~Ren, and J.~Sun}, {\em Delving deep into rectifiers:
  Surpassing human-level performance on imagenet classification}, in
  Proceedings of the IEEE international conference on computer vision, 2015,
  pp.~1026--1034.

\bibitem{heroux2005overview}
{\sc M.~A. Heroux, R.~A. Bartlett, V.~E. Howle, R.~J. Hoekstra, J.~J. Hu, T.~G.
  Kolda, R.~B. Lehoucq, K.~R. Long, R.~P. Pawlowski, E.~T. Phipps, et~al.},
  {\em An overview of the trilinos project}, ACM Transactions on Mathematical
  Software (TOMS), 31 (2005), pp.~397--423.

\bibitem{Khanjani2023}
{\sc Z.~Khanjani, G.~Watson, and V.~P. Janeja}, {\em Audio deepfakes: A
  survey}, Frontiers in Big Data, 5 (2023),
  \url{https://doi.org/10.3389/fdata.2022.1001063},
  \url{https://www.frontiersin.org/articles/10.3389/fdata.2022.1001063}.

\bibitem{Adam14}
{\sc D.~P. Kingma and J.~Ba}, {\em Adam: A method for stochastic optimization},
  2014, \url{https://doi.org/10.48550/ARXIV.1412.6980},
  \url{https://arxiv.org/abs/1412.6980}.

\bibitem{Adam2014}
{\sc D.~P. Kingma and J.~Ba}, {\em Adam: A method for stochastic optimization},
  arXiv preprint arXiv:1412.6980,  (2014).

\bibitem{luz2020}
{\sc I.~Luz, M.~Galun, H.~Maron, R.~Basri, and I.~Yavneh}, {\em Learning
  {{Algebraic Multigrid Using Graph Neural Networks}}}, Sept. 2020,
  \url{https://doi.org/10.48550/arXiv.2003.05744},
  \url{https://arxiv.org/abs/2003.05744}.

\bibitem{Moore2021}
{\sc N.~S. Moore, E.~C. Cyr, and C.~M. Siefert}, {\em Learning an algebriac
  multrigrid interpolation operator using a modified graphnet architecture},
  (2021), \url{https://doi.org/10.2172/1859673},
  \url{https://www.osti.gov/biblio/1859673}.

\bibitem{nocedal1999numerical}
{\sc J.~Nocedal and S.~J. Wright}, {\em Numerical optimization}, Springer,
  1999.

\bibitem{PyTorch2015}
{\sc A.~Paszke, S.~Gross, F.~Massa, A.~Lerer, J.~Bradbury, G.~Chanan,
  T.~Killeen, Z.~Lin, N.~Gimelshein, L.~Antiga, A.~Desmaison, A.~Kopf, E.~Yang,
  Z.~DeVito, M.~Raison, A.~Tejani, S.~Chilamkurthy, B.~Steiner, L.~Fang,
  J.~Bai, and S.~Chintala}, {\em Pytorch: An imperative style, high-performance
  deep learning library}, in Advances in Neural Information Processing Systems
  32, Curran Associates, Inc., 2019, pp.~8024--8035,
  \url{http://papers.neurips.cc/paper/9015-pytorch-an-imperative-style-high-performance-deep-learning-library.pdf}.

\bibitem{robbins1951stochastic}
{\sc H.~Robbins and S.~Monro}, {\em A stochastic approximation method}, The
  annals of mathematical statistics,  (1951), pp.~400--407.

\bibitem{saadBook}
{\sc Y.~Saad}, {\em Iterative Methods for Sparse Linear Systems}, Society for
  Industrial and Applied Mathematics, second~ed., 2003,
  \url{https://doi.org/10.1137/1.9780898718003},
  \url{https://epubs.siam.org/doi/abs/10.1137/1.9780898718003},
  \url{https://arxiv.org/abs/https://epubs.siam.org/doi/pdf/10.1137/1.9780898718003}.

\bibitem{scarselli_graph_2009}
{\sc F.~Scarselli, M.~Gori, A.~C. Tsoi, M.~Hagenbuchner, and G.~Monfardini},
  {\em The {Graph} {Neural} {Network} {Model}}, IEEE Transactions on Neural
  Networks, 20 (2009), pp.~61--80,
  \url{https://doi.org/10.1109/TNN.2008.2005605}.
\newblock Conference Name: IEEE Transactions on Neural Networks.

\bibitem{shukla2022scalable}
{\sc K.~Shukla, M.~Xu, N.~Trask, and G.~E. Karniadakis}, {\em Scalable
  algorithms for physics-informed neural and graph networks}, Data-Centric
  Engineering, 3 (2022), p.~e24.

\bibitem{taghibakhshi2022_OptimizationBasedAMG}
{\sc A.~Taghibakhshi, S.~MacLachlan, L.~Olson, and M.~West}, {\em
  Optimization-{{Based Algebraic Multigrid Coarsening Using Reinforcement
  Learning}}}, Jan. 2022, \url{https://arxiv.org/abs/2106.01854}.

\bibitem{taghibakhshietal2022_LearningInterfaceConditions}
{\sc A.~Taghibakhshi, N.~Nytko, T.~Zaman, S.~MacLachlan, L.~Olson, and
  M.~West}, {\em Learning {{Interface Conditions}} in {{Domain Decomposition
  Solvers}}}, May 2022, \url{https://arxiv.org/abs/2205.09833}.

\bibitem{taghibakhshi2023_MGGNN}
{\sc A.~Taghibakhshi, N.~Nytko, T.~U. Zaman, S.~MacLachlan, L.~Olson, and
  M.~West}, {\em {{MG-GNN}}: {{Multigrid Graph Neural Networks}} for {{Learning
  Multilevel Domain Decomposition Methods}}}, Jan. 2023,
  \url{https://doi.org/10.48550/arXiv.2301.11378},
  \url{https://arxiv.org/abs/2301.11378}.

\bibitem{tieleman2012lecture}
{\sc T.~Tieleman, G.~Hinton, et~al.}, {\em Lecture 6.5-rmsprop: Divide the
  gradient by a running average of its recent magnitude}, COURSERA: Neural
  networks for machine learning, 4 (2012), pp.~26--31.

\bibitem{wangetal2019_BackpropagationFriendlyEigendecomposition}
{\sc W.~Wang, Z.~Dang, Y.~Hu, P.~Fua, and M.~Salzmann}, {\em
  Backpropagation-{{Friendly Eigendecomposition}}}, June 2019,
  \url{https://doi.org/10.48550/arXiv.1906.09023},
  \url{https://arxiv.org/abs/1906.09023}.

\end{thebibliography}

\appendix
\section{GNN Architectural Details}

\subsection{GNN for Learning Jacobi Iteration}\label{app:arch:jacobi}
For the Section~\ref{sec:jacobi_nn} GNN, the architecture is
summarized in Table~\ref{tab:gnnArchJacobi}. As can be seen from the
table, there is a total of $1341$ trainable parameters.

\begin{table}[htpb]
    \centering
    \caption{GNN Architecture for Learning Jacobi Iteration}
    \label{tab:gnnArchJacobi}
    \begin{tabular}{cccccccc}
        GNN Layer & Update & Layer & num inputs & num outputs & bias & activation & num parameters \\
        \hline
        Layer 1 & $\phi _{v}$ & 1 & 5 & 50  & Yes & ReLU & $(5+1)\times 50$ \\
         & &              2 & 50 & 20 & Yes & ReLU & $(50+1)\times 20$ \\
         & &              3 & 20 & 1  & Yes &      & $(20+1)\times 1$
    \end{tabular}
\end{table}

\subsection{GNN for Diffusion Coefficients}\label{app:arch:diffusion}
The architecture of the best performing model among the 432 models
described in Section~\ref{sec:arch:diffusion} is summarized above in
Table~\ref{tab:gnnArchDiffCoeffs}. As described in
Section~\ref{sec:arch:diffusion}, the model consists of two GNN
layers; however the first of these layers as an encoder which performs
no aggregation.

The number of inputs in $\phi_{e}$ and $\phi_{v}$ may require some
explanation. From the algorithm, we see that $\phi_{e}$ takes as
arguments the attributes from the edge itself as well as its connected
vertices and the global variable. After the encoder, each of these is
a vector of length 32. Hence the input to $\phi_{e}$ in layer 1 is $32
\times 4 = 128$. Similarly, for the vertex update, $\phi_{v}$, the
aggregation function $\rho_{e\rightarrow v}$ returns four aggregated
attributes per input attribute. Hence, the input the the $\phi_{v}$
function includes 32 vertex attributes, $4 \times 32$ aggregated
attributes, and $32$ global attributes for a total of $192$
attributes.   

Finally, from the table we see that there are a total of $14002$
trainable parameters in the full network. 

\begin{table}[htpb]
    \centering
    \caption{GNN Architecture for Learning Diffusion Coefficients}
    \label{tab:gnnArchDiffCoeffs}
    \begin{tabular}{cccccccc}
        GNN Layer & Update & Layer & num inputs & num outputs & bias & activation & num parameters \\
        \hline
        Encoder & $\phi _{e}$ & 1 & 3  & 16 & Yes & ReLU & $(3+1)\times 16$ \\
                &             & 2 & 16 & 16 & Yes & ReLU & $(16+1)\times 16$ \\
                &             & 3 & 16 & 32 & Yes &      & $(16+1)\times 32$ \\
                & $\phi _{v}$ & 1 & 1  & 16 & Yes & ReLU & $(1+1)\times 16$ \\
                &             & 2 & 16 & 16 & Yes & ReLU & $(16+1)\times 16$ \\
                &             & 3 & 16 & 32 & Yes &      & $(16+1)\times 32$ \\
                & $\phi _{g}$ & 1 & 1  & 16 & Yes & ReLU & $(1+1)\times 16$ \\
                &             & 2 & 16 & 16 & Yes & ReLU & $(16+1)\times 16$ \\
                &             & 3 & 16 & 32 & Yes &      & $(16+1)\times 32$ \\
        Layer 1       & $\phi _{e}$ & 1 & 128& 32 & Yes & ReLU & $(128+1)\times 32$ \\
                &             & 2 & 32 & 32 & Yes &      & $(32+1)\times 32$ \\
                & $\phi _{v}$ & 1 & 192& 32 & Yes & ReLU & $(192+1)\times 32$ \\
                &             & 2 & 32 & 2  & Yes & Leaky ReLU & $(32+1)\times 2$
    \end{tabular}
\end{table}

\end{document}